\documentclass[10pt]{article}
\usepackage{amsmath, amsthm,  latexsym, amssymb, amsfonts, epsfig,
caption2, graphicx}
\usepackage[dvipdfm]{hyperref}

\def\cajita{\rule{5pt}{5pt}}
\renewcommand{\d}{\mathrm{d}}
\newcommand{\ep}{\varepsilon}
\newcommand{\ph}{\varphi}
\newcommand{\dw}{\downarrow}
\newcommand{\hy}{\widehat{Y}}
\newcommand{\hz}{\widehat{Z}}
\newcommand{\hu}{\widehat{u}}
\newcommand{\wt}{\widetilde}
\newcommand{\pr}{\prime}

\newcommand{\R}{\mathbb{R}}

\newcommand{\K}{\mathcal {K}}

\newtheorem{teorema}{Theorem}[section]
\newtheorem{proposicion}[teorema]
{Proposition}
\newtheorem{lema}[teorema]{Lemma}
\par
\newtheorem{definicion}[teorema]
{Definition}
\par
\newtheorem{corolario}[teorema]{Corollary}
\par
\par
\newtheorem{nota}[teorema]{Remark}
\par

\par
\newtheorem{Nota}[teorema]{Remark}
\par
\newtheorem{ejemplo}[teorema]
{Example}
\par

\par
\par

\begin{document}
\title{Semilinear Backward Doubly Stochastic Differential Equations and
SPDEs Driven by Fractional Brownian Motion with Hurst Parameter in
$(0,1/2)$}
\author{ Shuai Jing$^{a}$\thanks{Supported by the National Basic
Research  Program of China (973 Program) grant No. 2007CB814900
(Financial Risk) and Marie Curie Initial Training Network (ITN)
project: $"$Deterministic and Stochastic Controlled Systems and
Application$"$, FP7-PEOPLE-2007-1-1-ITN, No. 213841-2. E-mail:
\href{mailto:shuai.jing@univ-brest.fr}{shuai.jing@univ-brest.fr}} \
\ and Jorge A.~Le\'{o}n$^{b}$\thanks {Partially supported by the
CONACyT grant 98998. E-mail:
\href{mailto:jleon@ctrl.cinvestav.mx}{jleon@ctrl.cinvestav.mx}}\\
\small{$^a$ School of Mathematics, Shandong University, 250100,
Jinan, China}\\
\small{$^a$ D\'epartement de Math\'ematiques, Universit\'e de
Bretagne Occidentale,} \small{29285 Brest C\'edex, France }
 \\
\small{$^b$ Departamento de Control Autom\'atico, Cinvestav-IPN,
Apartado Postal 14--740, 07000 M\'{e}xico, D.F., Mexico}} \maketitle

\begin{abstract}
We study the existence of a unique solution to semilinear fractional
backward doubly  stochastic differential equation driven by a
Brownian motion and a fractional Brownian motion with Hurst
parameter less than 1/2. Here the stochastic integral with respect
to the fractional Brownian motion is the extended divergence
operator and the one with respect to Brownian motion is It\^o's
backward integral. For this we use the technique developed by
R.~Buckdahn \cite{bu} to analyze stochastic differential equations
on the Wiener space, which is based on the Girsanov theorem and the
Malliavin calculus, and we reduce the backward doubly stochastic
differential equation to a backward stochastic differential equation
driven by the Brownian motion. We also prove that the solution of
semilinear fractional backward doubly stochastic differential
equation defines the unique stochastic viscosity solution of a
semilinear stochastic partial differential equation driven by a
fractional Brownian motion.
\end{abstract}

\bigskip


\section{Introduction}
\label{sec:1}

This paper investigates semilinear fractional backward doubly
stochastic differential equations (BDSDEs) and semilinear stochastic
partial differential equations (SPDEs) driven by fractional Brownian
motion. Fractional Brownian motions (fBms) and backward stochastic
differential equations (BSDEs) have been extensively studied in
recent twenty years. However, up to now there are only few works
that combine both topics. Bender \cite{Be} considered a class of
linear fractional BSDEs and gave their explicit solutions. There are
two major obstacles depending on the properties of fBm: Firstly, the
fBm is not a semimartingale except for the case of Brownian motion
(Hurst parameter $H=1/2$), hence the classical It\^o calculus which
is based on semimartingales cannot be transposed directly to the
fractional case. Secondly, there is no martingale representation
theorem with respect to the fBm. However, such a martingale
representation property with respect to the Brownian motion is the
main tool in BSDE theory. Hu and Peng's paper \cite{HP} overcame the
second obstacle for the case of $H>1/2$ by using the
quasi-conditional expectation and by studying nonlinear fractional
BSDEs in a special case only.

Nevertheless, there are many papers considering stochastic
differential equations driven by fractional Brownian motion with
Hurst parameter $H>1/2$ (\cite{BHOZ}, \cite{Mi} and references
therein) or $H<1/2$ (\cite{Ls}), or covering both cases (\cite{Jm}).
For the case $H<1/2$, one of the main difficulties is how to
properly define the stochastic integral with respect to  the fBm. In
the paper of Cheridito and Nualart \cite{Cn}, and then generalized
by Le\'{o}n and Nualart \cite{Ln}, the authors have defined the
extended divergence operator which can be applied to the fBm for
$H<1/2$ as a special case. In this paper we will use such definition
for the stochastic integration with respect to the fBm, and then
apply the non-anticipating Girsanov transformation developed by
Buckdahn \cite{bu} to transform the semilinear fractional doubly
backward stochastic differential equation driven by the Brownian
motion $W$ and the fractional Brownian motion $B$:
\begin{equation}\label{1.1}
Y_t=\xi + \int_0^t f(s,Y_s,Z_s)\d s-\int_0^tZ_s\dw \d W_s + \int_0^t
\gamma_s Y_s \d B_s,\quad t\in\left[0,T\right],
\end{equation}
into a pathwise (in the sense of fBm) BSDE
\begin{equation}
\hy_t= \xi +\int^t_0 f\left(s,\hy_s \ep_s(T_s),
\hz_s\ep_s(T_s)\right)\ep_s^{-1}(T_s)\d s -\int_0^t \hz_s\dw \d
W_s,\quad t\in\left[0,T\right].
\end{equation}
More precisely, the solutions $(Y,Z)$ and $\left(\hy, \hz\right)$
are linked together by the following relations:
$$
\{(Y_t, Z_t), {t\in[0,T]}\} =\left \{\left(\hy_t(A_t)\ep_t,
\hz_t(A_t)\ep_t\right), t\in[0,T]\right\}
$$
and
$$
\left\{\left(\hy_t, \hz_t\right), {t\in[0,T]}\right\}
=\left\{\left(Y_t(T_t)\ep^{-1}_t(T_t),
Z_t(T_t)\ep^{-1}_t(T_t)\right), {t\in[0,T]}\right\},
$$
where $A_t$ and $T_t$ are Girsanov transformations. It is worth
noting that such kind of method was also used by Jien and Ma
\cite{Jm} to deal with fractional stochastic differential equations.

It is well known that the solution of a BSDE can be regarded as a
viscosity solution of an associated parabolic partial differential
equation (PDE) (cf. \cite{EPQ} \cite{Pa} and \cite{Pe}), and the
solution of BDSDE driven by two independent Brownian motions can be
regarded as a stochastic viscosity solution of an SPDE (cf.
\cite{BM} and \cite{PP2}). So it is natural to consider the
relationship between the solutions of our fractional BDSDE and the
associated SPDE. We show that the solution of the above fractional
BDSDE, which is a random field, is a stochastic viscosity solution
of an SPDE driven by our fractional Brownian motion. To be more
precise, the value function $u(t,x)$ defined by the solution of a
fractional BSDE over the time interval $[0,t]$ instead of $[0,T]$,
see (\ref{x-bdsde}), will be shown to possess a continuous version
and to be  the stochastic viscosity solution of the following
semilinear SPDE
\begin{equation}\label{1.3}
\left\{ \begin{array}{ll} \d u(t,x) = \left[ \mathcal{L}u(t,x) +
f\left(t,x,u(t,x),\nabla_x u(t,x)\sigma(x)\right)\right]\d t
+\gamma_t  u(t,x)\d B_t,& t\in[0,T],\\
u(0,x)=\Phi(x).&
\end{array}\right.
\end{equation}
Taking a Brownian motion instead of the fBm, equation (\ref{1.1})
becomes a classical BDSDE, which was first studied by Pardoux and
Peng \cite{PP2}. The associated stochastic viscosity solution of
SPDE (\ref{1.3}) (with $H=1/2$) was studied by Buckdahn and Ma
\cite{BM}. Let us point out that, unlike \cite{BM} considering the
stochastic integral with respect to $B$ ($H=1/2$) as the
Stratonovich one and using a Doss-Sussman transformation as main
tool, we have to do here with an extended divergence operator
($H<1/2$), which condemns us to use the Girsanov transformation as
main argument. However, this restricts us to semilinear equations.
We will investigate the general case in a forthcoming paper, but
with a different approach.

The paper is organized as follows: In Section \ref{sec:2} we recall
some preliminaries which will be used in what follows: Malliavin
calculus for fractional Brownian motion, the definition of extended
divergence operator and the Girsanov transformation. In Section
\ref{sec:3} we prove existence and uniqueness results for stochastic
differential equations driven by a fractional Brownian motion and
backward doubly stochastic differential equations driven by a
Brownian motion as well as a fractional Brownian motion. The
relationship between the stochastic viscosity solution of the
stochastic partial differential equation (\ref{1.3}) driven by
fractional Brownian motion and that of an associated pathwise
partial differential equation  is given in Section \ref{sec:4}.
\section{Preliminaries} \label{sec:2} \setcounter{equation}{0} The
purpose of this section is to describe the framework that will be
used in this paper. Namely, we introduce briefly the transformations
on the Wiener space, appearing in the construction of the solution
to our equations,  some preliminaries of the Malliavin calculus for
the fBm, and the left and right-sided fractional derivatives, which
are needed to understand the definition of the extension of the
divergence operator with respect to the fBm. Although the most of
results discussed in this section are known, we prefer to provide a
self--contained exposition for the convenience of the reader.
\subsection{Fractional calculus}
\label{sub:2.1}
For a detailed account on the fractional calculus theory, we refer,
for instance, to Samko et al. \cite{Sa}.

Let $T>0$ denote a positive time horizon, fixed throughout our
paper, and let $f:[0,T]\rightarrow \R$ be an integrable function,
and $\alpha \in (0,1)$. The {\it right--sided fractional integral}
of $f$ of order $\alpha$ is given by
$$
I^\alpha_{T-}(f)(x)=\frac{1}{\Gamma (\alpha)} \int^T_x
\frac{f(u)}{(u-x)^{1-\alpha}}\d u, \quad \hbox{\rm for a.a.}\;\;x\in
[0,T].
$$
Note that $I^\alpha_{T-}(f)$ is well-defined because the  Fubini
theorem implies that it is a function in $L^p ([0,T])$,  $p\geq 1$,
whenever $f\in L^p ([0,T])$.

We denote by $I^\alpha_{T-}(L^p),\  p\geq 1$, the family of all
functions $f\in L^p ([0,T])$ such that
\begin{equation}
\label{eq:2.2.1} f=I^\alpha_{T-}(\ph),
\end{equation}
for some $\ph \in L^p([0,T])$. Samko et al. \cite{Sa} (Theorem 13.2)
provide a characterization of the space $I^\alpha_{T-}(L^p)$, $p>1$.
Namely, a measurable function $f$ belongs to $I^\alpha_{T-}(L^p)$
(i.e., it satisfies $(\ref{eq:2.2.1}))$ if and only if $f\in L^p
([0,T])$ and the integral
\begin{equation}
\label{eq:2.2.2} \int^T_{s+\ep} \frac{f(s)-f(u)}{(u-s)^{1+\alpha}}\d
u
\end{equation}
converges in $L^p ([0,T])$ as $\ep \dw 0$. In this case a function
$\ph$ satisfying $(\ref{eq:2.2.1})$ coincides with the right--sided
fractional derivative
\begin{equation}
\label{eq:2.2.3} (D^\alpha_{T-}f)(s)=\frac{1}{\Gamma (1-\alpha)}
\biggl(\frac{f(s)}{(T-s)^\alpha}+\alpha \int^T_s
\frac{f(s)-f(u)}{(u-s)^{1+\alpha}}\d u\biggr),
\end{equation}
where the integral is the $L^p([0,T])$--limit of $(\ref{eq:2.2.2})$.
Moreover, it has also been shown in \cite{Sa} (Lemma 2.5) that there
is at most one solution $\ph$ to the equation $(\ref{eq:2.2.1})$.
Consequently, the inversion formulae
$$
I^\alpha_{T-}\left(D^\alpha_{T-}f\right)=f, \quad\hbox{\rm for all }
f\in I^\alpha_{T-}(L^p),
 $$
and
$$
D^\alpha_{T-}\left(I^\alpha_{T-}(f)\right)=f,\quad \hbox{\rm for all
} f\in L^1 ([0,T])$$ hold.

Similarly, the {\it left--sided fractional integral and the
derivative} of $f$ of order  $\alpha$, which are given,
respectively,  by
$$
I^\alpha_{0+}(f)(x)=\frac{1}{\Gamma (\alpha)} \int^x_0
\frac{f(u)}{(x-u)^{1-\alpha}}\d u, \quad \hbox{\rm for a.a.}\;\;x\in
[0,T],
$$
and
\begin{equation}\label{eq:lef-der}
(D^\alpha_{0+}f)(s)=\frac{1}{\Gamma (1-\alpha)}
\biggl(\frac{f(s)}{s^\alpha}+\alpha \int^s_0
\frac{f(s)-f(u)}{(s-u)^{1+\alpha}}\d u\biggr),
\end{equation}
satisfy the inversion formulae
$$
I^\alpha_{0+}\left(D^\alpha_{0+}f\right)=f, \quad\hbox{\rm for all }
f\in I^\alpha_{0+}(L^p),
 $$
and
$$
D^\alpha_{0+}\left(I^\alpha_{0+}(f)\right)=f,\quad \hbox{\rm for all
} f\in L^1 ([0,T]).
$$
\subsection{Fractional Brownian motion}
\label{sub:2.2} In this subsection we will recall some basic facts
of the fBm. The reader can consult Mishura \cite{Mi} and Nualart
\cite{Dn} and the  references therein for a more complete
presentation of this subject.

Henceforth $(\Omega,{\cal F}, P)$ and $W^0=\{ W^0_t: t\in[0,T]\}$
are the canonical Wiener space on the interval $[0,T]$ and  the
canonical Wiener process, respectively. This means, in particular,
that $\Omega=C_0([0,T])$ is the space of all continuous functions
$h: [0,T]\mapsto \R$ with $h(0)=0$, $\{W^0_t(\omega)=\omega(t),
t\in[0,T], \omega\in\Omega\}$ is the coordinate process on $\Omega$,
$P$ is the Wiener measure on $(\Omega, \mathcal{B}(\Omega))$ and
$\mathcal{F}$ is the completion of
$\mathcal{B}(\Omega)=\sigma\{W^0_s, s\in[0,T]\}$ with respect to
$P$. The noise under consideration is the process
$$B_t=\int_0^t K_H(t,s)\d W^0_s,\quad t\in[0,T],$$
where  $K_H$ is the kernel of the fBm with parameter $H\in(0,1/2)$.
That is,
$$
K_H(t,s)=C_H\left[\left(\frac ts\right)^{H-1/2}(t-s)^{H-1/2}
-(H-1/2)s^{1/2 - H}\int_s^t u^{H-3/2}(u-s)^{H-1/2}\d u\right],
$$
where $C_H =\sqrt{\frac{2H}{(1-2H)\beta (1-2H, H + {1}/{2})}}$. The
process  $B=\{B_t:t\in [0,T]\}$ is an fBm with Hurst parameter $H$,
defined on  $(\Omega, {\cal F}, P)$, i.e., $B$ is a Gaussian process
with zero mean and covariance function
$$
R_H(t,s):=E\left[B_tB_s\right]=\int_0^{s\wedge t}K_H(t,r)K_H(s,r)\d
r= \frac{1}{2}\left(t^{2H}+ s^{2H}-|t-s|^{2H}\right), \quad s, t\in
[0,T].
$$

Let ${\cal H}_H$ be the Hilbert space defined as the completion of
the space $L^e(0,T)$ of step functions on $[0,T]$ with respect to
the norm generated by the inner product
$$
\left\langle 1_{[0,t]}, 1_{[0,s]}\right\rangle_{{\cal
H}_H}=R_H(t,s)=E\left[B_tB_s\right],\quad t, s\in [0,T].
$$
From Pipiras and Taqqu \cite{Pt} (see also \cite{Dn}), it follows
that ${\cal H}_H$ coincides with the Hilbert space
$$
\Lambda^{1/2-H}_T
 :=\left\{f:[0,T]\rightarrow \R :
\exists \;\ph_f \in L^2 (0,T) \hbox{ such that }  f(u) =u^{1/2-H}
I^{1/2-H}_{T-}\left(s^{H-1/2} \ph_f (s)\right)(u)\right\}
$$
equipped with scalar product
$$
\left\langle f,g\right\rangle_{\Lambda_T^{1/2 -H}}=C_H^2
\Gamma(H+1/2)^2 \langle \ph_f, \ph_g\rangle_{L^2(0,T)}.
$$
So  $1_{[0,t]}\mapsto B_t$ can be extended to an isometry of
$\Lambda^{1/2-H}_{T}$ onto a Gaussian closed subspace of
$L^2(\Omega, {\cal F}, P)$. This isometry is denoted by $\ph\mapsto
B(\ph)$. Moreover, by the transfer principle (see Nualart
\cite{Dn}), the map ${\cal K}: \Lambda^{1/2-H}_T\rightarrow
L^2([0,T])$, defined by (see (\ref{eq:2.2.3}))
$$
({\cal K}\ph)(s)=C_H\Gamma(H+1/2)s^{1/2 -H}
\left(D_{T-}^{1/2-H}u^{H-1/2}\ph(u)\right)(s),\quad s\in[0,T],
$$
is an isometry such that
$$
B(\ph)=\int_0^T ({\cal K}\ph)(s)\d W^0_s\quad \hbox{\rm and}\quad
{\cal K}1_{[0,t]}=K_H(t,\cdot)1_{[0,t]},\ t\in[0,T].
$$

Using the properties of ${\cal K}$, Cheridito and Nualart \cite{Cn}
have extended the domain of the divergence operator with respect to
the fBm $B$. This extension of the divergence operator in the sense
of Malliavin calculus holds also true for some suitable Gaussian
processes (see Le\'on and Nualart \cite{Ln}). For $B$, this
extension is introduced as follows.

The following result  identifies the adjoint of the operator ${\cal
K}$ (see \cite{Ln}). It uses the left--sided fractional derivative
$D^\alpha_{0+}$ defined in $(\ref{eq:lef-der})$.

\begin{proposicion}
\label{def-adj} {\it Let $g:[0,t]\rightarrow \R$ be a function such
that $u\mapsto u^{1/2 -H} g(u)$ belongs to \\ $I^{1/2 -H}_{0+}(L^q
([0,a]))$ for some $q>(1/2 -H)^{-1}\vee H^{-1}$. Then, $g\in$ Dom
${\cal K}^*$, and for all $u\in [0,T]$},
$$
({\cal K}^*g)(u)=C_H\Gamma(H+ 1/2) u^{H-1/2} D^{1/2-H}_{0+}
\left(s^{1/2-H} g(s)\right)(u).
$$
\end{proposicion}

Let ${\cal S}$ (resp. ${\cal S}_{\cal K}$) denote the class of
smooth random variables of the form
\begin{equation}
\label{eq:fuc-sua} F=f(B(\ph_1), \ldots, B(\ph_n)),
\end{equation}
where $\ph_1,\ldots, \ph_n$ are in $\Lambda_T^{1/2 -H}$ (resp. in
the domain of the operator ${\cal K}^*{\cal K})$ and $f\in
C_p^\infty (\R^n)$. Here, $C^\infty_p (\R^n)$ is the set of
$C^\infty$ functions $f:\R^n \rightarrow \R$ such that $f$ and all
its partial derivatives have polynomial growth.

The derivative of the smooth random variable $F$ given by
$(\ref{eq:fuc-sua})$ is the $\Lambda_T^{1/2 -H}$--valued random
variable $DF$ defined by
$$
DF=\sum^n_{i=1}\frac{\partial f}{\partial x_i} (B(\ph_1),\ldots,
B(\ph_n))\ph_i .
$$

Now we can introduce the stochastic integral that we use in this
paper. It is an extended divergence operator with respect to $B$.

\begin{definicion}
\label{def-div} Let $u\in L^2 (\Omega , {\cal F}, P; L^2([0,T]))$.
We say that $u$ belongs to Dom $\delta$ if there exists $\delta (u)
\in L^2 (\Omega ) $ such that
\begin{equation}
\label{eq:rel-dua} E\left[\langle {\cal K}^* {\cal K}DF, u
\rangle_{L^2([0,T])}\right]=E\left[F\delta (u)\right], \quad
\hbox{\rm for every}\quad F\in {\cal S}_{\cal K}.
\end{equation}
In this case, the random variable $\delta (u)$ is called the
\hbox{\rm extended divergence} of $u$.
\end{definicion}
\begin{nota}
$i)$ In \cite{Ln} it is shown that the domain of ${\cal K}^*{\cal
K}$ is a dense subset of $\Lambda_T^{1/2 -H}$. Therefore, there is
at most one square integrable random variable $\delta (u)$ such that
$(\ref{eq:rel-dua})$ holds.

$(ii)$ In \cite{Cn} and \cite{Ln}, it is proven that the domain of
$\delta$ is bigger than that of the classical divergence operator,
which is defined by the chaos decomposition approach $($see Nualart
\cite{Nu}$)$.

$(iii)$ In Section $\ref{sec:3}$ and Section $\ref{sec:4}$, we use
the convention
$$
\int^t_0 u_s \d B_s=\delta (u 1_{[0,t]}),
$$
whenever $u1_{[0,t]}\in$ Dom $\delta$.

\end{nota}
\subsection{Girsanov transformations}
\label{sub:2.3}In this section we introduce the Girsanov
transformations on $\Omega$, which we consider in this paper.

In all which follows, we assume that $\gamma$ is a square-integrable
function satisfying the following hypothesis:
\begin{itemize}
\item[\textbf{(H1)}] $\gamma 1_{[0,t]}$ belongs to $\Lambda_T^{1/2 -H}$, for
every $t\in[0,T]$.
\end{itemize}

We emphasize that Le\'on and San Mart\'{\i}n \cite{Ls} (Lemma 2.3)
have shown the existence of square-integrable functions satisfying
the above hypothesis.

Now, for $t\in[0,T]$, we consider the  transformations on $\Omega$
of the form
$$
T_t(\omega)=\omega+\int_0^{{\cdot}\wedge t} ({\cal
K}\gamma1_{[0,t]})(r) \d r
$$
and
$$
A_t(\omega)=\omega-\int_0^{{\cdot}\wedge t} ({\cal
K}\gamma1_{[0,t]})(r) \d r.
$$

Notice that $A_tT_t$ and $T_tA_t$ are the identity operator of
$\Omega$, and that the Girsanov theorem leads to write
$$
B(\ph)(T_t)=B(\ph)+\int_0^T({\cal K}\gamma1_{[0,t]})(r)({\cal K}\ph)
(r)\d r=B(\ph)+\int_0^t\gamma_r({\cal K}^*{\cal K}\ph)(r)\d r,
$$
for all $\ph\in \hbox{\rm Dom }({\cal K}^*{\cal K})$, and
\[E[F]=E[F(A_t)\ep_t],\] with
$$
\ep_t=\exp\left(\int_0^t\gamma_r\d B_r -\frac12 \int_0^t(({\cal
K}\gamma1_{[0,t]})(r))^2 \d r\right).
$$

We will need the following estimate of the above exponential of the
integral with respect to the fractional Brownian motion:
\begin{lema}\label{lemma}
Let $\gamma: [0,T]\mapsto \R$, $\gamma\in L_p[0,T]\cap D^B_p[0,T]$,
for some $p> 1/H$, where $D_p^B[0,T]=\{\gamma:
[0,T]\mapsto\R|\int^T_0(\int^T_x\ph(x,t)\d t)^p \d x < \infty \}$
and set $\ph(x,t)=\frac{|\gamma(t)-\gamma(x)|}{(t-x)^{3/2-H}}
{1}_{\{0<x<t\leq T\}}$. Then there exists a constant $C(H,p)$ only
depending on $H$ and $p$, such that
\begin{equation}\label{est_eps}
E\left[ \exp\left\{\sup_{0\leq0\leq T}\left |\int^t_0 \gamma_s\d
B_s\right|\right\}\right]\leq 2\exp \left\{1/2\left(C(H,p)
G_p(0,T,\gamma)+4\sqrt{2}\right)^2\right\},
\end{equation}
where $G_p(0,T,\gamma): =\|\gamma\|_{L^p[0,T]}\cdot T^{H-1/p} +
T^{1/2- 1/p}\left(\int^T_0\left(\int^T_x\ph(x,t)\d t\right)^p
\right)^{1/p}.$
\end{lema}
\noindent\textit{Proof:} Let $I_T^{\ast}=\sup_{0\le t\le
T}|\int^t_0\gamma_s\d B_s|$. According to Lifshitz\cite{Lif} Theorem
1, P.141, and its corollary, for all $r>4\sqrt{2}D(T,
{\lambda}/{2})$, we have the inequality
\begin{equation}\label{est_I}
P\{I_T^{\ast}>r\}\leq 2\left(1-\Phi\left(\frac{r-4\sqrt{2}
D(T,{\lambda}/{2})}{\lambda}\right)\right),
\end{equation}
where
$\Phi(x)=\frac{1}{\sqrt{2\pi}}\int^x_{-\infty}\exp\{-{y^2}/{2}\}\d
y$, $D(T,{\lambda}/{2})$ is the Dudley integral (for more details,
we refer to Mishura\cite{Mi}), and
$\lambda^2=\sup_{t\in[0,T]}E\left[\left(\int^t_0 \gamma_s\d
B_s\right)^2\right]$.

Since $ E\left[\exp\{I_T^{\ast}\}\right]=1+\int^{+\infty}_0
\exp\{x\} P(I_T^{\ast}>x)\d x, $ by using the estimate of
(\ref{est_I}) we have
\begin{equation*}
\begin{aligned}
&E\left[\exp\{I_T^{\ast}\}\right]\\
= &1+\int^{4\sqrt{2}D(T,{\lambda}/{2})}_0 \exp\{x\}
P(I_T^{\ast}>x)\d x+\int_{4\sqrt{2}D(T,{\lambda}/{2})}^{+\infty}
\exp\{x\}P(I_T^{\ast}>x)\d x\\
\leq &\exp\left\{4\sqrt{2}D(T,{\lambda}/{2})\right\} +2
\int_{4\sqrt{2}D(T,{\lambda}/{2})}^{+\infty} \exp\{x\}\left
(1-\Phi\left(\frac{x-4\sqrt{2}D(T,{\lambda}/{2})}
{\lambda}\right)\right)\d x\\
=&\exp\left\{4\sqrt{2}D(T,{\lambda}/{2})\right\} +
2\int_{0}^{+\infty} \exp\left\{4\sqrt{2}D(T,{\lambda}/{2})+x\right\}
\left(1-\frac{1}{\sqrt{2\pi}}\int_{-\infty}^{ {x}/{\lambda}}
\exp\{-{y^2}/{2}\}\d y \right)\d x\\
\leq&2\exp\left\{{\lambda^2}/{2}+4\sqrt2 D(T,{\lambda}/{2})\right\}.
\end{aligned}
\end{equation*}
Moreover, from Theorem 1.10.6 of Mishura\cite{Mi} and its proof we
know that
$$
\lambda\leq C_1(H,p)G_p(0,T,\gamma)
$$
and
$$D(T,{\lambda}/{2})\leq
C_2(p) G_p(0,T,\gamma).
$$
By substituting them to the former inequality, we easily get the
wished result. \hfill$\cajita$
\section{Semilinear fractional SDEs and
fractional backward doubly SDEs} \label{sec:3}
\setcounter{equation}{0}

\subsection{Fractional anticipating semilinear equations}
\label{sub:3.1} In this subsection we discuss the existence and
uniqueness of solutions to anticipating semilinear equations  driven
by a fractional Brownian motion $B$ with Hurst parameter
$H\in(0,1/2)$. This type of equation was studied by Jien and Ma
\cite{Jm}, and since it motivates the approach in our work, we give
it in details.

We consider the fractional anticipating equation
\begin{equation}\label{eq:anteq}
X_t=\xi+\int_0^t b(s,X_s)\d s+\int_0^t\gamma_sX_s\d B_s,\quad
t\in[0,T].
\end{equation}
Here $\xi\in L^p(\Omega)$, $p>2$, and $b:\Omega\times[0,T]\times
\R\rightarrow \R$ is a measurable function such that:
\begin{itemize}
\item[\textbf{(H2)}] There exist  $\nu\in L^1([0,T])$, $\nu\ge 0$,
and a positive constant $L$ such that
\begin{eqnarray*}
|b(\omega,t,x)-b(\omega,t,y)|&\le &\nu_t|x-y|,\quad \int_0^T\nu_t\d t\le L,\\
|b(\omega,t,0)|&\le&L,
\end{eqnarray*}
for all $x,y\in\R$ and almost all $\omega\in\Omega$.
\end{itemize}

We observe that the above assumption guarantees that the pathwise
equation
$$
\zeta_t(x)=x+\int_0^t\ep^{-1}_s(T_s)b(T_s,s,\ep_s (T_s)\zeta_s(x))\d
s, \quad t\in[0,T],
$$
has a unique solution. Henceforth, we denote it by $\zeta$.

Now we can state the existence of a unique solution equation
(\ref{eq:anteq}):

\begin{teorema} Under Hypotheses \textbf{$($H1$)$} and \textbf{$($H2$)$},
the process
\begin{equation}\label{eq:sol-anteq}
X_t=\ep_t\zeta_t(A_t,\xi(A_t))
\end{equation}
is the unique solution in $L^2(\Omega\times[0,T])$ of the equation
$(\ref{eq:anteq})$, such that $\gamma X1_{[0,t]}\in$ Dom $\delta $,
for all $t\in[0,T].$
\end{teorema}

\noindent\textit{Proof:}  We first show that the process $X$ given
in (\ref{eq:sol-anteq}) is a solution of equation (\ref{eq:anteq}).
For this we first observe that $X$ belongs to
$L^2(\Omega\times[0,T])$ and we let  $F\in {\cal S}_{\cal K}$. Then,
by the integration by parts formula and the Girsanov theorem,
together with the fact that $\frac{\d F(T_t)}{\d t}= \gamma_t({\cal
K}^*{\cal K} DF)(T_t,t)$, we have
\begin{equation*}
\begin{aligned}
&E\left[FX_t-F\xi\right]
=E\left[F(T_t)\zeta_t(\xi)-F\zeta_0(\xi)\right]
=E\left[\int_0^t\frac{\d}{\d s}(F(T_s)\zeta_s(\xi))\d s\right]\\
=&\int_0^t\gamma_sE\left[({\cal K}^*{\cal K}
DF)(T_s,s)\zeta_s(\xi)\right]\d s +
\int_0^tE\left[F(T_s)\ep^{-1}_s(T_s)b(T_s,s,\ep_s
(T_s)\zeta_s(\xi))\right]\d s\\
=&\int_0^t\left(\gamma_sE\left[({\cal K}^*{\cal K}
DF)(s)X_s\right]+E\left[Fb(s,X_s)\right]\right)\d s.
\end{aligned}
\end{equation*}
Hence, since $\gamma X1_{[0,t]}$ is square integrable, Definition
\ref{def-div} implies that $\gamma X1_{[0,t]}$ belongs to Dom
$\delta$ and the equality in (\ref{eq:anteq}) holds, for all
$t\in[0,T]$.

Now we deal with the uniqueness of equation (\ref{eq:anteq}). For
this end, let $Y$ be another solution of equation (\ref{eq:anteq}),
$F\in {\cal S}_{\cal K}$ and $t\in[0,T]$. Then,
$$E\left[Y_tF(A_t)\right]=E\left[\xi F(A_t)\right]+E\left[
\int_0^tF(A_t)b(s,Y_s)\d s\right]+E\left[\int_0^t\gamma_sY_s
\left({\cal K}^*{\cal K}DF(A_t)\right)(s)\d s\right].$$ Therefore,
the  integration by parts formula, Fubini's theorem as well as the
fact that $\frac{\d F(A_s)}{\d s}=-\gamma_s({\cal K}^*{\cal
K}DF(A_s))(s)$ yield
\begin{equation*}
\begin{aligned}
&{E\left[Y_tF\right]} \\=&E\left[\xi F(A_t)\right]-E\left[\xi
\int_0^t\gamma_s
({\cal K}^*{\cal K}DF(A_s))(s)\d s\right]\\
&+E\left[\int_0^t F(A_s)b(s,Y_s)\d s\right]- E\left[\int_0^t
\gamma_r ({\cal K}^*{\cal K}DF(A_r))(r)\int_0^r
b(s,Y_s)\d s\d r\right]\\
&+E\left[\int_0^t\gamma_sY_s({\cal K}^*{\cal K}DF(A_s))(s)\d
s\right] -E\left[\int_0^t\int_0^r\gamma_r ({\cal K}^*{\cal
K}D(({\cal K}^*{\cal K}DF(A_r))(s)))(r)\gamma_sY_s\d s\d r \right].
\end{aligned}
\end{equation*}
Hence, by using that $Y$ is a solution of (\ref{eq:anteq}),
Definition \ref{def-div} and the relation
$$
({\cal K}^*{\cal K}D(({\cal K}^*{\cal K}DF(A_r))(s)))(r)= ({\cal
K}^*{\cal K}D(({\cal K}^*{\cal K}DF(A_r))(r)))(s),
$$
we obtain
$$
E\left[Y_tF(A_t)\right]=E\left[\xi F\right]+
E\left[\int_0^tF(A_s)b(s,Y_s)\d s\right].
$$
Consequently, by using the Girsanov theorem again, we get
$$
Y_t(T_t)\ep_t^{-1}(T_t) =\xi+\int_0^tb(T_s,s,Y_s(T_s)) \ep_s^{-1}
(T_s)\d s,
$$
which implies that $Y_t(T_t)\ep_t^{-1}(T_t)=\zeta(\xi).$ That is,
$Y_t=\ep_t\zeta_t(A_t,\xi(A_t))$, and therefore the proof is
complete.\hfill$\cajita$
\subsection{Fractional  backward doubly stochastic differential equations}
\label{sub:3.2}In this section we state some of the main results of
this paper. Namely, the existence and uniqueness of  backward doubly
stochastic differential equations driven by both a fractional
Brownian motion $B$ and a standard Brownian motion $W$.

Let $\{B_t: {0\leq t\leq T}\}$ be a one-dimensional fractional
Brownian motion with Hurst parameter $H\in(0,1/2)$, defined on the
classical Wiener space $(\Omega^{\pr},{\cal F} ^{B}, P^{B})$ with
$\Omega^\pr=C_0([0,T];\R)$, and
$\{W_t=(W^1_t,W^2_t,\cdots,W^d_t):{0\leq t\leq T}\}$  a
d-dimensional canonical Brownian motion defined on the classical
Wiener space $(\Omega^{\pr\pr},{\cal F} ^{W}, P^{W})$ with
$\Omega^{\pr\pr}=C_0([0,T]; \R^d)$. We put $(\Omega, {\cal F}^0,
P)=(\Omega^{\pr},{\cal F}^{B},P^{B}) \otimes(\Omega^{\pr\pr},{\cal
F} ^{W}, P^{W})$ and let $\mathcal {F}=\mathcal{F}^0\vee
\mathcal{N}$, where $\mathcal N$ is the class of the $P$-null sets.
We denote again by $B$ and $W$ their canonical extension from
$\Omega^{\pr}$ and $\Omega^{\pr\pr}$, respectively, to $\Omega$.

We let $\mathcal{F}_{t,T}^W = \sigma\{W_T-W_s, t\leq s\leq
T\}\vee\cal N$, $\mathcal{F}_t^{B}=\sigma\{B_s, 0\leq s\leq
t\}\vee\cal N,$ and
$\mathcal{G}_t=\mathcal{F}_{t,T}^W\vee\mathcal{F}_t^{B}, t\in[0,T].$
Let us point out that $\mathcal{F}_{t,T}^W$ is decreasing and
$\mathcal{F}_t^{B}$ is increasing in $t$, but $\mathcal{G}_t$ is
neither decreasing nor increasing. We denote the family of
$\sigma$-fields $\{\mathcal{G}_t\}_{0\leq t\leq T}$ by $\mathbb{G}$.
Moreover, we shall also introduce the backward filtrations
$\mathbb{H}=\{\mathcal {H}_t =\mathcal{F}_{t,T}^W \vee
\mathcal{F}_T^{B}\}_ {t\in[0,T]}$ and
$\mathbb{F}^W=\{\mathcal{F}_{t,T}^W\}_{t\in[0,T]}$.

Let ${\cal S}^\pr_{\cal K}$ denote the class of smooth random
variables of the form
\begin{equation}
\label{eq:fuc-sua2} F=f(B(\ph_1), \ldots, B(\ph_n), W(\psi_1),
\ldots, W(\psi_m)),
\end{equation}
where $\ph_1,\ldots, \ph_n$ are elements of the domain of the
operator ${\cal K}^*{\cal K}$, $\psi_1, \ldots, \psi_m\in
C([0,T],\R^d)$, $f\in C_p^\infty (\R^{n+m})$ and $n,m\geq1$. Here,
$C^\infty_p (\R^{n+m})$ is the set of all $C^\infty$ functions
$f:\R^{n+m} \rightarrow \R$ such that $f$ and all its partial
derivatives have polynomial growth.

The Malliavin derivative of the smooth random variable $F$  w.r.t.
$B$ is the $\Lambda_T^{1/2 -H}$--valued random variable $D^BF$
defined by
$$
D^B F=\sum^n_{i=1}\frac{\partial f}{\partial x_i} (B(\ph_1),\ldots,
B(\ph_n),W(\psi_1), \ldots, W(\psi_m))\ph_i,
$$
and the Malliavin derivative $D^W F$ of the smooth random variable
$F$ w.r.t. $W$ is given by
$$
D^W F=\sum^{m}_{i=1}\frac{\partial f}{\partial x_{n+i}} (B(\ph_1),
\ldots, B(\ph_n),W(\psi_1), \ldots, W(\psi_m))\psi_i.
$$
\begin{definicion}\label{def-div-ext}
$($Skorohod integral w.r.t.\,$B$. Extension of Definition
\ref{def-div}.$)$ We say that $u\in L^2 (\Omega \times[0,T])$
belongs to Dom $\delta^B$ if there exists a random variable
$\delta^B(u)\in L^2 (\Omega ) $ such that
\begin{equation*}
E\left[\langle {\cal K}^* {\cal K}D^BF, u
\rangle_{L^2(\left[0,T\right])}\right]=E\left[F\delta^B(u)\right],
\quad \hbox{\rm for all}\quad F\in {\cal S}^\pr_{\cal K}.
\end{equation*}
We call $\delta^B(u)$ the Skorohod integral with respect to $B$.
\end{definicion}
\begin{definicion}$($Skorohod integral w.r.t. $W$.$)$
We say that $u\in L^2 (\Omega \times[0,T])$ belongs to Dom
$\delta^W$ if there exists a random variable $\delta^W(u)\in L^2
(\Omega ) $ such that
\begin{equation*}
E\left[\int^T_0 (D^W_sF) u_s\d s\right]=E\left[F\delta^W(u)\right],
\quad \hbox{\rm for all}\quad F\in {\cal S}^\pr_{\cal K}.
\end{equation*}
We call $\delta^W(u)$ the Skorohod integral with respect to $W$.
\end{definicion}

For the the Skorohod integral with respect to $W$ we have, in
particular, the following well known result:
\begin{proposicion}\label{def-w}
Let  $u\in L^2 (\Omega \times[0,T])$ be $\mathbb{H}$-adapted. Then
the It\^o backward integral $\int^T_0 u_s\dw \d W_s$  coincides with
the Skorohod integral with respect to $W$: $$\int^T_0 u_s\dw \d W_s
=\delta^W(u).$$
\end{proposicion}

The extension to $\Omega$ of the operators $T_t$ and $A_t$
introduced in subsection \ref{sub:2.3} as acting over
$\Omega^{\pr}$, is done in a canonical way:
$$T_t(\omega^{\pr},\omega^{\pr\pr}):=(T_t\omega^{\pr},
\omega^{\pr\pr}),\quad
A_t(\omega^{\pr},\omega^{\pr\pr}):=(A_t\omega^{\pr},
\omega^{\pr\pr}), \quad \mathrm{for} \ (\omega^{\pr},
\omega^{\pr\pr})\in\Omega=\Omega^{\pr}\times\Omega^{\pr\pr}.$$

We denote by $L^2_{\mathbb{G}}(0,T;\R^n)$ (resp.,
$L^2_{\mathbb{H}}(0,T;\R^n)$) the set of $n$-dimensional measurable
random processes
$\{\ph_t,t\in[0,T]\}$ which satisfy:\\
i) $E\left[\int_{0}^T|\ph_t|^2\d t\right]<+\infty$, \\
ii) $\ph_t$ is $\mathcal{G}_t$- (resp., $\mathcal{H}_t$-)
measurable, for a.e. $t\in[0,T]$.

We also shall introduce a subspace of  $L^2_{\mathbb{G}}(0,T;\R^n)$,
which stems its importance from its invariance with respect to a
class of Girsanov transformations. Recalling the notation
$$
I^*_T:=\sup_{t\in[0,T]}\left|\int_0^t\gamma_sdB_s\right|
$$
from subsection \ref{sub:2.3}, we define
$L^{2,\ast}_{\mathbb{G}}(0,T;\R\times\R^d)$ to be the space of all
$\mathbb{G}$-adapted processes $(Y,Z)$  which are such that
$$
E\left[\exp\{p
I_T^{\ast}\}\int^T_0\left(|{Y}_t|^2+|{Z}_t|^2\right)\d t
\right]<\infty, \ for \ all \ p\ge 1.
$$
For the space $L^{2,\ast}_{\mathbb{G}}(0,T;\R\times\R^d)$ we have
the following invariance property:
\begin{proposicion}\label{prop_YZ}
For all processes $(Y,Z)\in L^{2,\ast}_\mathbb{G}(0,T;\R\times\R^d)$
we have:\\
i) $\left(\wt{Y}_t,\wt{Z}_t\right):=\left(Y_t(T_t)\ep_t^{-1}(T_t),
Z_t(T_t)\ep_t^{-1}(T_t)\right)\in
L^{2,\ast}_\mathbb{G}(0,T;\R\times\R^d)$
and\\
ii)
$\left(\overline{Y}_t,\overline{Z}_t\right):=\left(Y_t(A_t)\ep_t,
Z_t(A_t)\ep_t\right)\in L^{2,\ast}_\mathbb{G}(0,T;\R\times\R^d).$
\end{proposicion}
\noindent\textit{Proof:} Since the proofs of i) and ii) are similar,
we only prove i):

For the case of $(\wt{Y},\wt{Z})$, from the Girsanov transformation
and Lemma \ref{lemma} we have
\begin{equation*}
\begin{aligned}
&E\left[\exp\{p
I^{\ast}_T\}\int^T_0\left(\left|\wt{Y}_t\right|^2+\left|\wt{Z}_t\right|^2\right)\d
t\right] \\
=&\int^T_0E\left[\exp\{p I^{\ast}_T\}(|{Y}_t(T_t)|^2+|{Z}_t(T_t)|^2)
\ep_t^{-2}(T_t)\right]\d
t\\
=&\int^T_0E\left[\exp\{p I^{\ast}_T(A_t)\}(|{Y}_t|^2+|{Z}_t|^2)
\ep_t^{-1}\right]\d t\\
\leq& E\left[\exp\left\{p \sup_{0\leq t\leq T
}\left|\int^t_0\gamma_s\d B_s\right|+{\sup_{{0\le t\le T}\atop{0\leq
r \le T} }}\int^T_0(\mathcal{K}\gamma 1_{[0,t]})(s)
(\mathcal{K}\gamma1_{[0,r]}) (s)\d s
\right\}\int^T_0(|{Y}_t|^2+|{Z}_t|^2)
\ep_t^{-1}\d t\right]\\
\leq& C E\left[\exp\{(p+1) I^{\ast}_T\}\int^T_0(|{Y}_t|^2+|{Z}_t|^2)
\d
t\right]\\
<&+\infty,\ \ for \ all\  p\ge1.
\end{aligned}
\end{equation*}
Hence, the proof is complete.\hfill$\cajita$

\bigskip
We now consider the following type of backward doubly stochastic
differential equation driven by the Brownian motion $W$ and the
fractional Brownian motion $B$:
\begin{equation}\label{BDSDE}
Y_t=\xi + \int_0^t f(s,Y_s,Z_s)\d s-\int_0^tZ_s\dw \d W_s + \int_0^t
\gamma_s Y_s \d B_s,\quad t\in\left[0,T\right].
\end{equation}
Here $\xi\in L^2(\Omega,\mathcal{F}_{0,T}^W, P) $ and
$f:\Omega^{\pr\pr}\times\left[0,T\right]\times\R\times\R^d\mapsto\R$
is a measurable function such that
\begin{itemize}
\item[\textbf{(H3)}] \textbf{i)}. $f(\cdot,t, y, z)$ is
$\mathcal{F}_{t,T}^W$-measurable, for all $t\in [0, T]$, and for all
$(y, z)\in \R \times \R^d$;
\item[\textbf{(H3)}] \textbf{ii)}. $f(\cdot,0,0)\in L^2(\Omega\times
\left[0,T\right]);$
\item[\textbf{(H3)}] \textbf{iii)}. There exists a constant $C\in \R^+$ such
that for all $(y_1,z_1),(y_2,z_2)\in\R\times\R^d$,
$$|f(t,y_1,z_1)-f(t,y_2,z_2)|\leq C(|y_1-y_2|+|z_1-z_2|),\quad
a.e., a.s.$$
\end{itemize}
\begin{nota}
Let us refer to some special cases of the above BDSDE:

$i)$ If $\gamma=0$, equation (\ref{BDSDE}) becomes a classical BSDE
$($Pardoux and Peng \cite{PP3}$)$ with a unique solution $(Y,Z)\in
L^2_{\mathbb{F}^W}$ $ (0,T; \R\times \R^d);$

$ii)$ If $\xi\in \R$ and $f:[0,T]\times \R\times \R^d\rightarrow \R$
are deterministic, we can choose $Z=0$ and $Y\in L^2(\Omega'\times
[0,T])$ with $\gamma Y1_{[0,t]}\in\mbox{Dom}(\delta^B),\,
t\in[0,T],$ as the unique solution of the fractional SDE
$$ Y_t=\xi+\int_0^tf(s,Y_s,0)\d s+\int_0^t\gamma_sY_s\d B_s,\, \,
t\in[0,T],$$ which can be solved due to subsection \ref{sub:3.1}.

$iii)$ Pardoux and  Peng \cite{PP2} considered  backward doubly SDEs
in the case that $B$ and $W$ are two independent Brownian motions
and in a nonlinear framework.

\end{nota}

We let $\widetilde{\Omega}^\pr:=
\{\omega^\pr\in\Omega^\pr|I_T^*(\omega^\pr)<\infty\}$, which
satisfies $P^B(\widetilde{\Omega}^\pr)=1$ from Lemma \ref{lemma}. We
first establish the following theorem:
\begin{teorema}
\label{thm-pwBSDE} For all $\omega^\pr\in \widetilde{\Omega}^\pr$,
the backward stochastic differential equation
\begin{equation}\label{bsde}
\hy_t(\omega^\pr,\cdot)= \xi +\int^t_0
f\left(s,\hy_s(\omega^\pr,\cdot) \ep_s(T_s,\omega^\pr),
\hz_s(\omega^\pr,\cdot)
\ep_s(T_s,\omega^\pr)\right)\ep_s^{-1}(T_s,\omega^\pr)\d s -\int_0^t
\hz_s(\omega^\pr,\cdot)\dw \d W_s,
\end{equation}
$t\in[0,T]$, has a unique solution $\left(\hy(\omega^\pr,\cdot),
\hz(\omega^\pr,\cdot)\right) \in
L^{2}_{\mathbb{F}^W}(0,T;\R\times\R^d).$

Moreover, putting
$\left(\hy_t(\omega^\pr,\cdot),\hz_t(\omega^\pr,\cdot)\right):
=(0,0),$ for $\omega^\pr\in\widetilde{\Omega}^{\pr c}$, the random
variable $\Big(\hy_t(\omega^\pr,\omega^{\pr\pr})$, $
\hz_t(\omega^\pr,\omega^{\pr\pr})\Big)$ is jointly measurable in
$(\omega^\pr,\omega^{\pr\pr})$, and $\left(\hy,\hz\right)\in
L^{2,*}_{\mathbb{G}}(0,T;\R\times\R^d)$.

Furthermore, there exists a positive constant  $C$ $($only depending
on the $L^2$-norm of $\xi$ and $\K \gamma 1_{[0,t]}$, $L^2$-bound of
$f(\cdot,0,0)$ and the Lipschitz constant of $f)$ such that, for all
$\omega^\pr\in\widetilde{\Omega}^\pr$:
\begin{equation}\label{estimation}
E^W\left[\sup_{t\in[0,T]}\left|\hy_t(\omega^\pr,\cdot)\right|^2+
\int^T_0\left|\hz_t(\omega^\pr,\cdot)\right|^2\d t\right]\le
C\exp\{2I^{\ast}_T(\omega^\pr)\}.
\end{equation}
\end{teorema}
\noindent\textit{Proof:} We put  $F_s(\omega^\pr,y,z)=
f\left(s,y\ep_s(T_s,\omega^\pr),z\ep_s(T_s,\omega^\pr)\right)
\ep_s^{-1}(T_s,\omega^\pr)$, $s\in[0,T], (y,z)\in\R\times\R^d$,
$\omega^\pr\in\widetilde{\Omega}^\pr$. Obviously, \\
i) $F_s(\cdot,y,z)$ is $\mathcal{G}_s$-measurable and
$F_s(\omega^\pr,\cdot,y,z)$ is $\mathcal{F}_{s,T}^W$-measurable,
$\omega^\pr\in\widetilde{\Omega}^\pr$;\\
ii) $F_s(\omega^\pr,\omega^{\pr\pr},y,z)$ is Lipschitz  in $(y,z)$,
uniformly
with respect to $(s,\omega^\pr,\omega^{\pr\pr})$;\\
iii) $|F_s(\omega^\pr,\cdot, 0,0)|\le
C|f(\cdot,s,0,0)|\exp\{I_T^*(\omega^\pr)\}$.

Using  $F_s$, equation (\ref{bsde}) can be rewritten as follows:
\begin{equation}\label{bsde1}
\hy_t(\omega^\pr)= \xi +\int^t_0
F_s\left(\omega^\pr,\hy_s(\omega^\pr,\cdot),\hz_s(\omega^\pr,\cdot)\right)\d
s -\int_0^t\hz_s(\omega^\pr,\cdot)\dw \d W_s, \
t\in\left[0,T\right],\ \omega^\pr\in\widetilde{\Omega}^\pr.
\end{equation}

{\bf Step 1}: We begin by proving the {\it existence}: From the
conditions i)-iii) and standard BSDE arguments (see: Pardoux and
Peng \cite{PP3}) we know that, for all
$\omega^\pr\in\widetilde{\Omega}^\pr$, there is a unique solution
$\left(\hy(\omega^\pr,\cdot),\hz(\omega^\pr,\cdot)\right)\in
L^2_{\mathbb{F}^W}(0,T;\R\times\R^d)$. On the other hand, the joint
measurability of $F_s$ with respect to
$(\omega^\pr,\omega^{\pr\pr})$ allows to show that, extended to
$\Omega^\pr\times\Omega^{\pr\pr}$ by putting
$\left(\hy_t(\omega^\pr,\cdot),\hz_t(\omega^\pr,\cdot)\right):
=(0,0), \ \omega^\pr\in\widetilde{\Omega}^{\pr c}$, the process
$\left(\hy,\hz\right)$ is $\mathbb{H}$-adapted. Let us show that
$\left(\hy,\hz\right)\in L^2_{\mathbb{H}}(0,T;\R\times\R^d)$. For
this end, it suffices to prove (\ref{estimation}).

Let $\omega^\pr\in\widetilde{\Omega}^{\pr}$ be arbitrarily fixed. By
applying It\^o's formula to $\left|\hy_t\right|^2$ we have at
$\omega^\pr$, $P^W$-a.s.,
\begin{equation*}
\begin{aligned}
\d \left|\hy_t\right|^2=2\hy_t\left(F_t\left(\hy_t,\hz_t\right)\d t
- \hz_t\dw\d W_t\right)-\left|\hz_t\right|^2\d t.
\end{aligned}
\end{equation*}
It follows that at $\omega^\pr$, $P^W$-a.s.,
\begin{equation*}
\begin{aligned}
\left|\hy_t\right|^2+\int^t_0\left|\hz_s\right|^2\d
s&=\xi^2+2\int^t_0\hy_sF_s\left(\hy_s,\hz_s\right)\d
s - 2\int^t_0\hy_s\hz_s\dw\d W_s\\
&\leq\xi^2+\int^t_0\left(C_1\left|\hy_s\right|^2+\frac{1}{2}\left|\hz_s\right|^2+
C_2|F_s(0,0)|^2\right)\d s - 2\int^t_0\hy_s\hz_s\dw\d W_s.
\end{aligned}
\end{equation*}
Hence, at $\omega^\pr$, $P^W$-a.s.,
\begin{equation*}
\begin{aligned}
\left|\hy_t\right|^2+\frac12\int^t_0\left|\hz_s\right|^2\d s
\leq\xi^2+\int^t_0\left(C_1\left|\hy_s\right|^2+
C_3\ep_s^{-2}\right)\d s - 2\int^t_0\hy_s\hz_s\dw\d W_s.
\end{aligned}
\end{equation*}
Taking the expectation with respect to $P^W$, we notice that
$$
\begin{aligned}
&E^W\left[\left(\int_0^T\left|\hy_t(\omega^\pr,\cdot)
\hz_t(\omega^\pr,\cdot)\right|^2\d t\right)^{1/2}\right] \\
\le& \left(E^W\left[\sup_{t\in[0,T]}
\left|\hy_t(\omega^\pr,\cdot)\right|^2\right]\right)^{1/2}
\left(E^W\left[\int_0^T \left|\hz_t(\omega^\pr,\cdot)\right|^2\d
t\right]\right)^{1/2}<+\infty.
\end{aligned}
$$
Consequently, $E^W\left[\int^t_0\hy_s(\omega^\pr,\cdot)
\hz_s(\omega^\pr,\cdot)\dw \d W_s\right]=0,$ and by taking the
conditional expectation with respect to $\mathcal{F}_T^B$, we obtain
\begin{equation*}
\begin{aligned}
E^W\left[\left|\hy_t(\omega^\pr,\cdot)\right|^2
+\frac12\int^t_0\left|\hz_s(\omega^\pr,\cdot)\right|^2\d s\right]
\leq E^W\left[\xi^2\right] +
\int^t_0C_1E^W\left[\left|\hy_s(\omega^\pr,\cdot)\right|^2\right]\d
s+ C_4\exp\{2I^{\ast}_T(\omega^\pr)\}.
\end{aligned}
\end{equation*}
Thus, from Gronwall's inequality, we have
\begin{equation*}
\begin{aligned}
E^W\left[\left|\hy_t(\omega^\pr,\cdot)\right|^2\right]\leq
\left(E\left[\xi^2\right]+
C_4\exp\left\{2I^{\ast}_T\left(\omega^\pr\right)\right\}\right)\exp\{C_1t\},
t\in[0,T],
\end{aligned}
\end{equation*}
which, combined with the previous estimate, yields
$$
E^W\left[\int^T_0\left(\left|\hy_t(\omega^\pr,\cdot)\right|^2
+\left|\hz_t(\omega^\pr,\cdot)\right|^2\right)\d t\right]\le
C\exp\{2I^{\ast}_T(\omega^\pr)\}.
$$
In order to get the estimate (\ref{estimation}), it suffices now to
estimate
$E^W\left[\sup_{t\in[0,T]}\left|\hy_t(\omega^\pr,\cdot)\right|^2\right]$
by using equation (\ref{bsde1}), Burkh\"{o}lder-Davis-Gundy's
inequality and the above estimate.

To prove that $\left(\hy,\hz\right)$ belongs even to
$L^2_{\mathbb{G}}(0,T;\R\times\R^d)$, we have to prove the {\it
uniqueness} in $L^2_{\mathbb{H}}(0,T;\R\times\R^d)$.

{\bf{Step 2:}} We suppose that $\left(\hy^1,\hz^1\right)$ and
$\left(\hy^2,\hz^2\right)$ are two solutions of equation
(\ref{bsde1}) belonging to  $ L^2_{\mathbb{H}}(0,T;\R\times\R^d)$
(Notice that $L^2_{\mathbb{G}}(0,T;\R\times\R^d)\subset
L^2_{\mathbb{H}}(0,T;\R\times\R^d))$. Putting $\Delta
\hy_t=\hy_t^1-\hy_t^2$ and $\Delta \hz_t=\hz_t^1-\hz_t^2$, we have
\begin{equation*}
\Delta \hy_t = \int^t_0
\left[F_s\left(\hy_s^1,\hz_s^1\right)-F_s\left(\hy_s^2,
\hz_s^2\right)\right] \d s -\int_0^t \Delta \hz_s\dw \d W_s, \quad
t\in[0,T].
\end{equation*}
By applying It\^o's formula to $\left|\Delta \hy_t\right|^2$, we get
that
\begin{equation*}
\begin{aligned}
E\left|\Delta\hy_t\right|^2+E\left[\int^t_0
\left|\Delta\hz_s\right|^2\d s\right]&=
2E\left[\int^t_0\Delta\hy_s\left[F_s\left(\hy_s^1,
\hz_s^1\right)-F_s\left(\hy_s^2,\hz_s^2\right)\right]
\d s\right]\\
&\leq E\left[\int^t_0\left[2(C_0+C_0^2)\left|\Delta\hy_s\right|^2 +
\frac{1}{2}\left|\Delta\hz_s\right|^2\right]\d s\right],
\end{aligned}
\end{equation*}
and finally from Gronwall's lemma, we conclude that $\Delta\hy_t=0$,
$\Delta\hz_t=0,$ a.s., a.e.

{\bf{Step 3:}} Let us now show that $\left(\hy,\hz\right)$ is not
only in $L^2_{\mathbb{H}}(0,T;\R\times\R^d)$ but even in
$L^2_{\mathbb{G}}(0,T;\R\times\R^d)$. For this we consider for an
arbitrarily given $\tau\in\left[0,T\right]$ equation (\ref{bsde1})
over the time interval $\left[0,\tau\right]$:
\begin{equation}\label{bsde2}
\hy_t^\tau= \xi +\int^t_0 F_s\left(\hy_s^\tau,\hz_s^\tau\right)\d s
-\int_0^t\hz_s^\tau\dw \d W_s, \quad t\in\left[0,\tau\right].
\end{equation}
Let
$\mathcal{H}_t^{\tau}:=\mathcal{F}_{t,T}^W\vee\mathcal{F}_{\tau}^{B},$
$t\in[0,\tau]$. Then
$\mathbb{H}^{\tau}=\{\mathcal{H}_t^{\tau}\}_{t\in[0,\tau]}$ is a
filtration with respect to which $W$ has the martingale
representation property. Since $F_t(y,z)$ is $\mathcal{G}_t$- and
hence also $\mathcal{H}_t^{\tau}$-measurable, $\d t$ a.e. on
$[0,\tau]$, it follows from the classical BSDE theory that BSDE
(\ref{bsde2}) admits a  solution $\left(\hy,\hz\right) \in
L^2_{\mathbb{H}^{\tau}}(0,\tau;\R\times\R^d).$   Due to the first
step this solution is unique in $
L^2_{\mathbb{H}^{\tau}}(0,T;\R\times\R^d).$ Hence,
$\left(\hy_t,\hz_t\right)=\left(\hy^{\tau}_t,\hz^{\tau}_t\right)$,
$\d t$ a.e., for $t<\tau$. Consequently, $\left(\hy_t,\hz_t\right)$
is $\mathcal{H}_t^\tau$- measurable, $\d t$ a.e., for $t<\tau.$
Therefore, letting $\tau\dw t$ we can deduce from the right
continuity of the filtration $\mathbb{F}^B$ that
$\left(\hy,\hz\right)$ is $\mathbb{G}$-adapted.

It still remains to prove that $\left(\hy,\hz\right)\in
L^{2,*}_{\mathbb{G}}(0,T;\R\times\R^d)$.

{\bf Step 4:}   For the proof that $\left(\hy,\hz\right)$ belongs to
$L^{2,\ast}_\mathbb{G}(0,T;\R\times\R^d)$ we notice that, by the
above estimates and Lemma \ref{lemma},
\begin{equation*}
\begin{aligned}
&E\left[\exp\{p
I^{\ast}_T\}\int^T_0\left(\left|\hy_t\right|^2+\left|\hz_t\right|^2\right)\d
t\right] = E\left[\exp\{p I^{\ast}_T\}E\left[
\int^T_0\left(\left|\hy_t\right|^2+\left|\hz_t\right|^2\right)\d
t|\mathcal{F}^B_T\right]\right]\\
\leq&E\left[C\exp\{(2+p)I^{\ast}_T\}\right]<\infty, \ \ for \ all\ \
p\ge1.
\end{aligned}
\end{equation*}
Hence, the proof is complete.\hfill$\cajita$

\begin{corolario}
The process $\left(\hy,\hz\right)$ given by Theorem \ref{thm-pwBSDE}
is the unique solution of equation $(\ref{bsde})$ in
$L^{2,*}_{\mathbb{G}}(0,T; \R\times\R^d)$.
\end{corolario}

Now we state the main result of this subsection:
\begin{teorema}\label{thm-BDSDE}
$1)$ Let $\left(\hy,\hz\right)\in
L^{2,\ast}_\mathbb{G}(0,T;\R\times\R^d)$ be a solution of BSDE
$(\ref{bsde})$. Then
$$
\left\{(Y_t, Z_t), {t\in[0,T]}\right\} =
\left\{\left(\hy_t(A_t)\ep_t,
\hz_t(A_t)\ep_t\right),{t\in[0,T]}\right\}\in
L^{2,\ast}_\mathbb{G}(0,T;\R\times\R^d)
$$
is a solution of equation $(\ref{BDSDE})$ with $\gamma Y
1_{[0,t]}\in Dom \ \delta^{B},$ for all $t\in[0,T]$.

$2)$ Conversely, given an arbitrary solution $(Y,Z)\in
L^{2,\ast}_\mathbb{G}(0,T;\R\times\R^d)$ of the equation
$(\ref{BDSDE})$ with $\gamma Y 1_{[0,t]}\in Dom \ \delta^{B}$, for
all $t\in[0,T]$,  the process
$$
\left\{\left(\hy_t, \hz_t\right), {t\in[0,T]}\right\}
=\left\{\left(Y_t(T_t)\ep^{-1}_t(T_t),
Z_t(T_t)\ep^{-1}_t(T_t)\right), {t\in[0,T]}\right\}\in
L^{2,\ast}_\mathbb{G}(0,T;\R\times\R^d)
$$
is a solution of BSDE $(\ref{bsde})$.
\end{teorema}

From the Theorems \ref{thm-pwBSDE} and \ref{thm-BDSDE}  we can
immediately conclude the following

\begin{corolario}
The solution of equation $(\ref{BDSDE})$ in
$L^{2,\ast}_{\mathbb{G}}(0,T;\R\times\R^d)$ exists and is unique.
\end{corolario}

\noindent\textit{Proof $($of Theorem \ref{thm-BDSDE}$)$:}  We first
prove that, given the solution $\left(\hy,\hz\right)$ of equation
(\ref{bsde}), $(Y,Z)$ defined in the theorem solves (\ref{BDSDE}).
For this we notice that for $F$ being an arbitrary but fixed element
of $\mathcal{S}^{\pr}_{\mathcal{K}}$, from Girsanov transformation
and from the equation (\ref{bsde}), it follows
\begin{equation*}
\begin{aligned}
&E\left[FY_t-F\xi\right]=E \left[F(T_t)\hy_t-F\hy_0\right]\\=
&E\left[F(T_t)\hy_0+ F(T_t)\int^t_0F_s\left(\hy_s,\hz_s\right)\d s -
F(T_t)\int^t_0\hz_s\dw\d W_s- F\hy_0\right].
\end{aligned}
\end{equation*}
We recall that $E\left[F(T_t)\int^t_0\hz_s\dw\d
W_s\right]=E\left[\int^t_0D_s^W F(T_t)\hz_s\d s\right]$. Thus, from
the fact that $\frac{\d}{\d t} F(T_t)=\gamma_t
(\mathcal{K}^*\mathcal{K}D^B F)$ $(t,T_t)$, we have
\begin{equation*}
\begin{aligned}
E\left[FY_t-F\xi\right]
=&E\left[\hy_0\int^t_0\gamma_s(\mathcal{K}^*\mathcal{K}D^B
F)(s,T_s)\d s \right]\\
&+E\left[\int^t_0 F(T_s)F_s\left(\hy_s, \hz_s\right)\d s
+\int^t_0\int^t_s \gamma_r(\mathcal{K}^*\mathcal{K}D^B F)(r,T_r)\d r
F_s\left(\hy_s,\hz_s\right)\d s\right]\\
&-E\left[\int^t_0 D^W_s F(T_s)\hz_s\d s +\int^t_0\int^t_s D^W_s
\gamma_r(\mathcal{K}^*\mathcal{K}D^B F)(r,T_r)\d r\hz_s\d s\right].
\end{aligned}
\end{equation*}
Moreover, from Fubini's theorem, the definition of the Skorohod
integral with respect to $W$, and from Proposition \ref{def-w}, we
obtain that
\begin{equation*}
\begin{aligned}
&E\left[\int^t_0\int^t_s D^W_s \gamma_r(\mathcal{K}^*\mathcal{K}D^B
F)(r,T_r)\d r\hz_s\d s\right] =\int^t_0\gamma_rE\left[\int^r_0
D_s^{W} (\mathcal{K}^*\mathcal{K}D^B F)(r,T_r)\hz_s \d s\right]\d r\\
=&\int^t_0\gamma_rE\left[(\mathcal{K}^*\mathcal{K}D^B
F)(r,T_r)\int^r_0\hz_s\dw\d W_s\right]\d r.
\end{aligned}
\end{equation*}
Thus, by applying  the inverse Girsanov transformation as well as
Fubini's theorem, we obtain
\begin{equation*}
\begin{aligned}
E\left[FY_t-F\xi\right] =&E\left[\int^t_0 F(T_s)F_s\left(\hy_s,
\hz_s\right)\d
s  - \int^t_0 D^W_s F(T_s) \hz_s\d s\right]\\
&+E\left[\int^t_0\gamma_s(\mathcal{K}^*\mathcal{K}D^B
F)(s,T_s)\left(\hy_0+\int^s_0 F_r\left(\hy_r,\hz_r\right)\d r -
\int^s_0
\hz_r\dw\d W_r\right)\d s\right] \\
=&E\left[\int^t_0 Ff(s,{Y}_s, {Z}_s)\d s\right] - E\left[\int^t_0
D^W_s F {Z}_s\d s \right] + E\left[ \int^t_0 \gamma_s
(\mathcal{K}^*\mathcal{K}D^B F)(s){Y}_s\d s \right],
\end{aligned}
\end{equation*}
where, for the latter expression, we have used that $\hy_s$ is a
solution of $(\ref{bsde})$. Since $Z\in L^{2,\ast}_ \mathbb{G}
(0,T;\R\times\R^d) \left(\subset L^{2}_\mathbb{G}
(0,T;\R\times\R^d)\right)$, it holds
$$
E\left[\int^t_0 D^W_s F {Z}_s\d s\right]=E\left[F \int^t_0
{Z}_s\dw\d W_s\right].
$$
Consequently, $E\left[\int^t_0(\mathcal{K}^*\mathcal{K}D^B
F)(s)\gamma_s{Y}_s\d s \right]= E \left[ F \left\{ Y_t - \xi -
\int^t_0 f(s,Y_s,Z_s)\d s+\int^t_0 Z_s\dw\d W_s\right\}\right]$
holds for all $F\in \mathcal{S}_{\mathcal{K}}^\pr.$ From Proposition
\ref{prop_YZ} we know that both, $\left(\hy,\hz\right)$ and $(Y,Z)$,
belong to $L^{2,\ast}_\mathbb{G}(0,T;\R\times\R^d).$ Consequently,
$Y_t-\xi-\int^t_0f(s,Y_s,Z_s)\d s+\int^t_0 Z_s\dw\d W_s\in
L^2(\Omega,{\cal F},P)$, for all $t\in[0,T]$. Moreover, using
(\ref{estimation}),
\[\begin{aligned}
&E\left[\int^T_0\left|\gamma_r1_{[0,t]}(r)Y_r\right|^2\d r\right] =
E\left[\int^T_0\left|\gamma_r1_{[0,t]}(r)\hy_r(A_r)\right|^2\ep_r^2\d
r\right]\\
=&
E\left[\int^T_0\left|\gamma_r1_{[0,t]}(r)\hy_r\right|^2\ep_r(T_r)\d
r\right] \le
CE\left[\int^T_0\left|\gamma_r1_{[0,t]}(r)\hy_r\right|^2\exp\{I_T^*\}\d
r\right]\\
\le&
CE\left[\int^T_0\left|\gamma_r1_{[0,t]}(r)\right|^2\sup_{r\in[0,T]}
\left|\hy_r\right|^2 \exp\{I_T^*\}\d r\right] \le
C\int^T_0\left|\gamma_r1_{[0,t]}(r)\right|^2
E\left[\exp\{CI_T^*\}\right]\d r\\
\le & C\int^T_0\left|\gamma_r1_{[0,t]}(r)\right|^2\d r < \infty.
\end{aligned}
\]
Thus, according to the definition of the Skorohod integral with
respect to $B$ we then conclude $\gamma Y 1_{[0,t]}\in Dom \
\delta^{B}$ and
$$
\int^t _0 \gamma_s Y_s \d B_s =Y_t-\left(\xi + \int^t_0 f(s,Y_s,Z_s)
\d s - \int_0^t Z_s\dw \d W_s\right), t\in[0,T].
$$
Hence $Y_t=\xi + \int ^t_0f(s,Y_s,Z_s)\d s - \int_0^t Z_s\dw \d
W_s+\int^t _0 \gamma_s Y_s \d B_s , t\in[0,T].$

\medskip

Let us prove now the second assertion of the Theorem. For this end
we let $(Y,Z)\in L^{2,\ast}_{\mathbb{G}}(0,T;\R\times\R^d)$ be an
arbitrary solution of equation (\ref{BDSDE}) and $F$ be an arbitrary
but fixed element of $\mathcal{S}^{\pr}_{\mathcal{K}}$. Then we have
\begin{equation*}
\begin{aligned}
&E \left[Y_tF(A_t)\right]\\=&E\left[\xi
F(A_t)\right]+E\left[\int_0^t F(A_t)f(s,Y_s,Z_s)\d s\right]
-E\left[F(A_t)\int_0^t Z_s \dw \d W_s \right]+ E
\left[F(A_t)\int^t_0 \gamma_s Y_s
\d B_s \right]\\
=&I_1+I_2-I_3+I_4,
\end{aligned}
\end{equation*}
where, using the fact that $\frac{\d}{\d t}F(A_t) = -\gamma_t
(\mathcal{K}^* \mathcal{K}D^B F(A_t))(t)$ and Fubini's theorem,
\begin{equation*}
\begin{aligned}
I_1=&E\left[\xi F\right]-E\left[\xi \int^t_0\gamma_s
(\mathcal{K}^\ast\K D^BF(A_s))(s)\d s\right],
\end{aligned}
\end{equation*}
\begin{equation*}
\begin{aligned}
I_2=&E\left[\int_0^t F(A_s)f(s,Y_s,Z_s)\d s\right] -E\left[ \int_0^t
\gamma_r (\K^\ast\K D^B F(A_r))(r)\int^r_0f(s,Y_s,Z_s)\d s\d
r\right].
\end{aligned}
\end{equation*}
From the duality between the It\^o backward integral and the
Malliavin derivative $D^W$ (recall that $Z$ is square integrable),
we get
\begin{equation*}
\begin{aligned}
I_3=
&E \left[\int_0^tZ_s D_s^W F(A_s)\d s
\right]-E\left[\int^t_0Z_sD^W_s\left[\int_s^t\gamma_r(\K^\ast\K
D^BF(A_r))(r)\d r\right]\d s\right]\\
=&E \left[\int_0^t Z_s D_s^W F(A_s)\d s
\right]-E\left[\int^t_0\int_0^r Z_s\gamma_rD^W_s(\K^\ast\K
D^BF(A_r))(r)\d s\d r\right]\\
=&E \left[\int_0^t Z_s D_s^W F(A_s)\d s
\right]-E\left[\int^t_0\gamma_s(\K^\ast\K
D^BF(A_s))(s)\int_0^s Z_r\dw\d W_r\d s\right].\\
\end{aligned}
\end{equation*}
Moreover, from the duality between the integral w.r.t. $B$ and $D^B$
(observe that  $\gamma Y 1_{[0,t]}\in Dom\ \delta^B$) we obtain
\begin{equation*}
\begin{aligned}
I_4
=&E\left[\int^t_0 \gamma_s Y_s (\K^\ast\K D^BF(A_s))(s)\d s\right]-
E\left[\int^t_0 \int^r_0 \gamma_r(\K^\ast\K D^B(\K^\ast\K
D^BF(A_r)(r))(s))\gamma_s Y_s\d s\d r
\right]\\
= &E\left[\int^t_0 \gamma_s Y_s (\K^\ast\K D^BF(A_s))(s)\d s\right]-
E\left[\int^t_0 \gamma_s (\K^\ast\K D^BF(A_s))(s)\int^s_0 \gamma_r
Y_r\d B_r\d s \right].
\end{aligned}
\end{equation*}
Consequently, using the fact that $(Y,Z)$ is a solution of equation
(\ref{BDSDE}), i.e.,
$$
Y_t=\xi + \int_0^t f(s,Y_s,Z_s)\d s -\int_0^tZ_s\dw \d W_s +
\int_0^t \gamma_s Y_s \d B_s, \quad t\in[0,T],
$$
we obtain that $ E \left[Y_tF(A_t)\right]=E\left[\xi F\right] +
E\left[\int_0^t F(A_s)f(s,Y_s,Z_s)\d s\right] - E\left[\int_0^t Z_s
D_s^W F(A_s)\d s \right]. $ Therefore, by applying Girsanov
transformation again and taking into account the arbitrariness of
$F\in\mathcal{S}_{\mathcal{K}}^{\pr}$, it follows that for all
$F\in\mathcal{S}_{\mathcal{K}}^{\pr},$
\begin{equation*}
\begin{aligned}
&E\left[\int^t_0 Z_s(T_s)\ep_s^{-1}(T_s) D_s^W F\d s\right]
=E\left[F\left\{\xi-Y_t(T_t)\ep_t^{-1}(T_t)
+\int^t_0f(s,Y_s(T_s),Z_s(T_s))\ep_s^{-1}(T_s)\d s\right\}\right].
\end{aligned}
\end{equation*}
Now, since according to Proposition \ref{prop_YZ},
$(Y_t(T_t)\ep_t^{-1}(T_t),Z_t(T_t)\ep_t^{-1}(T_t))\in
L^{2,\ast}_{\mathbb{G}}(0,T;\R\times\R^d)$, we have\\
$Y_t(T_t)\ep_t^{-1}(T_t)- \xi - \int^t_0 f(s,Y_s(T_s),Z_s(T_s))
\ep_s^{-1}(T_s) \d s \in L^{2}(\Omega,{\cal G}_t,P)$. Therefore,
\begin{equation}\label{bsde3}
\begin{aligned}
Y_t(T_t)\ep_t^{-1}(T_t) =\xi + \int_0^t
f\left(s,Y_s\left(T_s\right),
Z_s\left(T_s\right)\right)\ep_s^{-1}(T_s)\d s - \int^t_0
Z_s(T_s)\ep^{-1}_s(T_s)\dw \d W_s,  a.s.,
\end{aligned}
\end{equation}
for all $t\in[0,T]$, which means that
$(\hy,\hz):=\left\{Y_t(T_t)\ep^{-1}_t(T_t),
Z_t(T_t)\ep^{-1}_t(T_t),t\in[0,T]\right\}$ is a solution of equation
(\ref{bsde}). Hence, the proof is complete. \hfill$\cajita$
\section{The associated stochastic partial differential equations}
\label{sec:4} \setcounter{equation}{0}

In this section we will use the following standard assumptions:
\begin{itemize}
\item[\textbf{(H4)}] \textbf{i)}. The functions $\sigma:\R^d\to\R^{d\times d},$
$b:\R^d\to\R^d$ and $\Phi:\R^d\to\R$ are  Lipschitz.
\item[\textbf{(H4)}] \textbf{ii)}. The function
$f:\left[0,T\right]\times\R^d\times\R\times\R^d\mapsto\R$ is
continuous, $f(t,\cdot,\cdot,\cdot)$ is Lipschitz, uniformly with
respect to $t$ and $f(\cdot,0,0,0)\in L^2(\Omega\times
\left[0,T\right])$.
\end{itemize}

We denote by $(X_s^{t,x})_{0\leq s \leq t}$ the unique solution of
the following stochastic differential equation:
\begin{equation}
\label{x-sde} \left\{
\begin{array}{l}
\d X_s^{t,x}=- b\left(X_s^{t,x}\right)\d
s -\sigma\left(X_s^{t,x}\right)\dw\d  W_s , \quad s\in\left[0,t\right],\\
X_t^{t,x}=x\in\R^d.
\end{array}\right.
\end{equation}
We note that this equation looks like a backward stochastic
differential equation, but due to the backward It\^o integral, the
SDE (\ref{x-sde}) is indeed a classical forward SDE. Under our
standard assumptions on $\sigma$ and $b$, it has a unique strong
solution which is $\{\mathcal{F}^W_{s,t}\}_{s\in[0,t]}$-adapted. The
following result provides some standard estimates for the solution
of equation (\ref{x-sde}) (cf. {\cite{PP1}} Lemma 2.7).
\begin{lema}\label{lemma_sol_x}
Let ${X^{t,x}}=\{X_s^{t,x},s\in[0,t]\}$ be the solution of the SDE
$(\ref{x-sde})$. Then

$(i)$ There exists a continuous version of $X^{t,x}$ such that
$(s,x)\mapsto X_s^{t,x}$ is locally H\"older
$\big(C^{\alpha,2\alpha}$, for all $\alpha\in (0,1/2)\big)$;

$(ii)$ For all $q\geq1$, there exists $C_q>0$ such that, for $t,
t^{\pr}\in[0,T]$ and $x, x^{\pr}\in\R^d$,
\begin{equation}
E\left[\sup_{0\leq s\leq t}\left|X_s^{t,x}\right|^q\right]\leq C_q
(1+|x|^q),
\end{equation}
\begin{equation}
E\left[\sup_{0\leq s\leq T}\left|X_{s\wedge t}^{t,x}-X_{s\wedge
t^\pr}^{t^{\pr},x^{\pr}} \right|^q\right]\leq C_q
\left[(1+|x|^q+|x^{\pr}|^q)|t-t^{\pr}|^{q/2}+|x-x^{\pr}|^q\right].
\end{equation}

\end{lema}

With the forward SDE we associate a BDSDE with driving coefficient
$f$:
\begin{equation}
\label{x-bdsde} Y_s^{t,x}=\Phi\left(X_0^{t,x}\right)
+\int^s_0f\left(r,X_r^{t,x},Y_r^{t,x},Z_r^{t,x}\right)\d r -
\int^s_0Z_r^{t,x}\dw\d  W_r +\int^s_0\gamma_rY_r^{t,x}\d B_r, \quad
s\in\left[0,t\right].
\end{equation}
By Theorem \ref{thm-pwBSDE} and Theorem \ref{thm-BDSDE}, the above
BDSDE has a unique solution $\left(Y^{t,x},Z^{t,x}\right)$ given by
$$\left(Y_s^{t,x}, Z_s^{t,x}\right) = \left(\hy_s^{t,x}(A_s)\ep_s,
\hz_s^{t,x}(A_s)\ep_s\right),  s\in[0,t],$$ where, for all
$\omega^\pr\in \widetilde{\Omega}^\pr$, $P^W$-a.s.,
\begin{equation}\label{x-bsde}
\begin{aligned}
\hy_s^{t,x}\left(\omega^\pr,\cdot\right)= \Phi\left(X_0^{t,x}\right)
+&\int^s_0
f\left(r,X_r^{t,x},\hy_r^{t,x}\left(\omega^\pr,\cdot\right)
\ep_r(T_r,\omega^\pr), \hz_r^{t,x}\left(\omega^\pr,\cdot\right)
\ep_r(T_r,\omega^\pr)\right)\ep_r^{-1}(T_r,\omega^\pr)\d r\\
&-\int_0^s \hz_r^{t,x}\left(\omega^\pr,\cdot\right)\dw \d W_r,\quad
s\in[0,t].
\end{aligned}
\end{equation}

Pardoux and Peng \cite{PP3} \cite{PP2} have studied BSDEs in a
Markovian context in which the driver $F_r(x,y,z)$ is deterministic;
here, in our framework the driver
$$
F_r(\omega^\pr,x,y,z):=f\left(r,x,y\ep_r(T_r,\omega^\pr),
z\ep_r(T_r,\omega^\pr)\right) \ep_r^{-1}(T_r,\omega^\pr), \
(r,x,y,z)\in[0,T]\times\R^d\times\R\times\R^d,
$$
is random but it depends only on $B$ and is independent of the
driving Brownian motion $W$. In the following we shall define
$X^{t,x}_s$, $\hy^{t,x}_s$ and $\hz^{t,x}_s$ for all
$(s,t)\in[0,T]^2$ by setting $X^{t,x}_s=x$,
$\hy^{t,x}_s=\hy^{t,x}_t$ and $\hz^{t,x}_{s}=0$, for $t<s\leq T$.

\bigskip
We have the following standard estimates for the solution:
\begin{lema}\label{lemma_est}
For  all $p\geq1$, there exists a constant $C_p\in \R_+$ such that
for all $(t,x),$ $(t',x')\in[0,T]\times \R^d,$
$\omega^\pr\in\widetilde{\Omega}^\pr$, $P^W$-a.s.,
\begin{equation}\label{est:sup_Y}
E^W\left[\sup_{0\leq s \leq
t}\left|\hy_s^{t,x}\left(\omega^\pr,\cdot\right)\right|^p + \left(
\int^t_0 \left|\hz_s^{t,x}\left(\omega^\pr,\cdot\right)\right|^2 \d
s \right)^{{p}/{2}}\right]\leq C_p\left(1+|x|^p\right)
\exp\left\{pI^*_T(\omega^\pr)\right\};
\end{equation}
\begin{equation}\label{est:Yt-Yt'}
\begin{aligned}
&E^W\left[\sup_{0\leq s\leq t\wedge t^\pr }
\left|\hy_s^{t,x}\left(\omega^\pr,\cdot\right)-
\hy_s^{t^\pr,x^\pr}\left(\omega^\pr,\cdot\right)\right|^p
+\left(\int^{ t\wedge
t^\pr}_0\left|\hz_s^{t,x}\left(\omega^\pr,\cdot\right)-
\hz_s^{t^\pr,x^\pr}\left(\omega^\pr,\cdot\right)\right|^2\d
s\right)^{p/2}\right] \\
\leq& C_p\exp\left\{pI^*_T(\omega^\pr)\right\}
\left(|x-x^\pr|^p+(1+|x|^p+|x^\pr|^p)|t-t^\pr|^{p/2}\right).
\end{aligned}
\end{equation}
\end{lema}
\noindent\textit{Proof}:  Let us fix any
$\omega^\pr\in\widetilde{\Omega}^\pr$. For $p\geq 2$, by applying
It\^o's formula to $|\hy_s^{t,x}|^p$, it follows that,  $P^W$-a.s.,
for $s\in[0,t]$,
$$
\begin{aligned}
&\left|\hy_s^{t,x}\left(\omega^\pr,\cdot\right)\right|^p+\frac{p(p-1)}{2}\int^s_0
\left|\hy_r^{t,x}\left(\omega^\pr,\cdot\right)\right|^{p-2}
\left|\hz_r^{t,x}\left(\omega^\pr,\cdot\right)\right|^2 \d r\\
=&\left|\Phi\left(X^{t,x}_0\right)\right|^p\\
&+p\int^s_0\left|\hy_r^{t,x}\left(\omega^\pr,\cdot\right)\right|^{p-2}
\hy_r^{t,x}\left(\omega^\pr,\cdot\right)\left(F_r\left(\omega^\pr,X_r^{t,x},
\hy_r^{t,x}\left(\omega^\pr,\cdot\right),
\hz_r^{t,x}\left(\omega^\pr,\cdot\right)\right)\d
r+\hz_r^{t,x}\left(\omega^\pr,\cdot\right)\dw\d W_r\right).
\end{aligned}
$$
Let $0\leq s\leq t^\pr\le t\le T$. We take the conditional
expectation with respect to $\mathcal{F}_{t^\pr,t}^W$ on both sides
of the above equality, and we obtain
$$
\begin{aligned}
&E^W\left[\left|\hy_s^{t,x}\left(\omega^\pr,\cdot\right)\right|^p
|\mathcal{F}_{t^\pr,t}^W\right] +E^W\left[\frac{p(p-1)}{2}\int^s_0
\left|\hy_r^{t,x}\left(\omega^\pr,\cdot\right)\right|^{p-2}
\left|\hz_r^{t,x}\left(\omega^\pr,\cdot\right)\right|^2 \d
r|\mathcal{F}_{t^\pr,t}^W\right]\\
=&E^W\left[\left|\Phi\left(X^{t,x}_0\right)\right|^p|\mathcal{F}_{t^\pr,t}^W\right]+
p E^W \bigg[ \int^s_0
\left|\hy_r^{t,x}\left(\omega^\pr,\cdot\right)\right|^{p-2}
\hy_r^{t,x}\left(\omega^\pr,\cdot\right)\\
&\quad\quad\quad\quad\quad\quad\times F_r\left(\omega^\pr,
X_r^{t,x},\hy_r^{t,x}\left(\omega^\pr,\cdot\right),\hz_r^{t,x}
\left(\omega^\pr,\cdot\right)\right)
\d r|\mathcal{F}_{t^\pr,t}^W\bigg]\\
\leq &E^W\left[\left|\Phi\left(X^{t,x}_0\right)\right|^p
|\mathcal{F}_{t^\pr,t}^W\right]+ p E^W \bigg[ \int^s_0 \Big( C_p
\left|\hy_r^{t,x}\left(\omega^\pr,\cdot\right)\right|^p
+C_p\left|X_r^{t,x}\right|^p \\&\quad\quad\quad\quad\quad+
\frac{p(p-1)} {4}
\left|\hy_r^{t,x}\left(\omega^\pr,\cdot\right)\right|^{p-2}
\left|\hz_r^{t,x}\left(\omega^\pr,\cdot\right)\right|^2+
C_p\ep_r^p(\omega^\pr)\Big)\d r|\mathcal{F}_{t^\pr,t}^W\bigg].
\end{aligned}
$$
Thus, from Gronwall's inequality and \textbf{(H4)} we have
$$
\left|\hy_s^{t,x}\left(\omega^\pr,\cdot\right)\right|^p\le
C_p\left(E^W\left[\sup_{0\leq r\leq
s}\left|X^{t,x}_r\right|^p|\mathcal{F}_{s,t}^W\right]
+\exp\left\{pI^*_T(\omega^\pr)\right\}\right)\leq
C_p\left(1+\left|X^{t,x}_s\right|^p\right)\exp\left\{pI^*_T(\omega^\pr)\right\}.
$$
Consequently, by using Doob's inequality, we get from the
arbitrariness of $p\ge 1$:
\[
E^W\left[\sup_{0\le s\le
t}\left|\hy^{t,x}_s\left(\omega^\pr,\cdot\right)\right|^p\right]\le
C_p\left(1+|x|^p\right)\exp\left\{pI_T^*(\omega^\pr)\right\}.
\]
The first result follows from Burkh\"older-Davis-Gundy inequality
applied to
$\left|\int_0^{t^\pr}\hz_s^{t,x}\left(\omega^\pr,\cdot\right)\dw \d
W_s\right|^p$ (see, e.g. \cite{PP1}).

Concerning the second assertion, without loss of generality, we can
suppose $t\geq t^{\pr}$. Let $0\leq s \le t^{\pr\pr}\leq t^{\pr}$.

Using an argument similar to that developed by Pardoux and Peng
\cite{PP1}, we see that, for some constants $\theta>0$ and $C>0$,
\begin{equation*}
\begin{aligned}
&E^W\left[\left|\hy_s^{t,x}\left(\omega^\pr,\cdot\right)
-\hy_s^{t^\pr,x^\pr}\left(\omega^\pr,\cdot\right)
\right|^p|\mathcal{F}_{t^{\pr\pr},t}^W\right]\\
&+ \theta
E^W\left[\int^s_0\left|\hy_r^{t,x}\left(\omega^\pr,\cdot\right)-
\hy_r^{t^\pr,x^\pr}\left(\omega^\pr,\cdot\right)\right|^{p-2}
\left|\hz_r^{t,x}\left(\omega^\pr,\cdot\right)
-\hz_r^{t^\pr,x^\pr}\left(\omega^\pr,\cdot\right)\right|^2\d
r|\mathcal{F}_{t^{\pr\pr},t}^W\right]\\
\leq& CE^W\left[\left|\Phi\left(X_0^{t,x}\right)
-\Phi\left(X_0^{t^\pr,x^\pr}\right)\right|^p|\mathcal{F}_{t^{\pr\pr},t}^W\right]\\
&+CE^W\left[\int^s_0\left[\ep^{-p}_r\left(T_r,\omega^\pr\right)
\left|X_r^{t,x}-X_r^{t^\pr,x^\pr}\right|^p+
\left|\hy_r^{t,x}\left(\omega^\pr,\cdot\right)
-\hy_r^{t^\pr,x^\pr}\left(\omega^\pr,\cdot\right) \right|^p\right]\d
r|\mathcal{F}_{t^{\pr\pr},t}^W\right]\\
\leq& CE^W\left[\left|\Phi\left(X_0^{t,x}\right)
-\Phi\left(X_0^{t^\pr,x^\pr}\right)\right|^p|\mathcal{F}_{t^{\pr\pr},t}^W\right]
+C\exp\left\{pI^*_T\left(\omega^\pr\right)\right\}\left(E\left[\int^s_0
\left|X_r^{t,x}-X_r^{t^\pr,x^\pr}\right|^{2p}\d
r|\mathcal{F}_{t^{\pr\pr},t}^W \right]\right)^{1/2}\\
&+CE^W\left[\int^s_0 \left|\hy_r^{t,x}\left(\omega^\pr,\cdot\right)
-\hy_r^{t^\pr,x^\pr}\left(\omega^\pr,\cdot\right)\right|^p\d
r|\mathcal{F}_{t^{\pr\pr},t}^W\right].
\end{aligned}
\end{equation*}
Consequently, from Gronwall's lemma and according to Lemma
\ref{lemma_sol_x},
\[
\left|\hy^{t,x}_s\left(\omega^\pr,\cdot\right)-\hy^{t,x}_s\left(\omega^\pr,\cdot\right)\right|^p\leq
C_p
\left|X^{t,x}_s-X^{t^\pr,x^\pr}_s\right|^p\exp\left\{pI_T^*\left(\omega^\pr\right)\right\},
0\le s\le t^\pr\le t,
\]
and
\[
E^W\left[\sup_{0\le s\le
t^\pr}\left|\hy^{t,x}_s\left(\omega^\pr,\cdot\right)
-\hy^{t,x}_s\left(\omega^\pr,\cdot\right)\right|^p\right] \le C_p
\exp\left\{pI_T^*\left(\omega^\pr\right)\right\}\left(\left(1+|x|^p+|x^\pr|^p\right)
|t-t^\pr|^{p/2}+|x-x^\pr|^p\right).
\]
Finally, with the help of the Burkh\"older-Davis-Gundy inequality
together with Lemma \ref{lemma}, we deduce that for all $p\geq 2$,
there exists $C_p$ such that
\begin{equation*}
\begin{aligned}
&E^W\left[\sup_{0\leq s\leq  t^\pr }
\left|\hy_s^{t,x}\left(\omega^\pr,\cdot\right)
-\hy_s^{t^\pr,x^\pr}\left(\omega^\pr,\cdot\right)\right|^p
+\left(\int^{ t^\pr}_0\left|\hz_s^{t,x}\left(\omega^\pr,\cdot\right)
-\hz_s^{t^\pr,x^\pr}\left(\omega^\pr,\cdot\right)\right|^2\d
s\right)
^{p/2}\right] \\
\leq& C_p\exp\left\{pI^*_T\left(\omega^\pr\right)\right\}
\left(|x-x^\pr|^p+(1+|x|^p+|x^\pr|^p)|t-t^\pr|^{p/2}\right).
\end{aligned}
\end{equation*}
The case $p\geq1$ follows easily from the case $p\geq 2$. This
completes the proof of the lemma.\hfill$\cajita$

\bigskip
We now introduce the random field: $\hu(t,x)=\hy_t^{t,x},
(t,x)\in[0,T]\times\R^d$, which has the following regularity
properties:
\begin{proposicion}\label{prop_meas}
The random field $\hu(t,x)$ is ${\cal F}^{B}_t$-measurable and  we
have $\hy_s^{t,x}\left(\omega^\pr,\cdot\right)$$=
\hy_s^{s,X_s^{t,x}}\left(\omega^\pr,\cdot\right)$
$=\hu(\omega^\pr,s,X_s^{t,x}),$ $P^W$-a.s.,  $0\leq s\leq t\leq T,$
$\omega^\pr\in\widetilde{\Omega}^\pr$.
\end{proposicion}
\noindent\textit{Proof:} From Theorem \ref{thm-pwBSDE} with terminal
time $t$, we know that $\hy_s^{t,x}$ is
$\mathcal{F}_{s,t}^W\vee\mathcal{F}_s^B$-measurable. Hence
$\hu(t,x)=\hy^{t,x}_t$ is
$\mathcal{F}_{t,t}^W\vee\mathcal{F}_t^B$-measurable. By applying
Blumenthal zero-one law we deduce that $\hu$ is
$\mathcal{F}_t^B$-measurable and independent of $W$. The second
assertion is a direct result from the uniqueness property of the
solutions of (\ref{x-sde}) and (\ref{x-bsde})(cf. \cite{EPQ}).
\hfill$\cajita$
\begin{lema}\label{continuity}
The process $\{\hy_s^{t,x};(s,t)\in[0,T]^2,x\in\R^d\}$ possesses a
continuous version. Moreover, $|\hu(t,x)|\leq
C\exp\{I_T^{\ast}\}(1+|x|),$ $P$-a.s.
\end{lema}
\noindent\textit{Proof:} Recall that, for $s\in[0,t],$

$$\widehat{Y}^{t,x}_s=\Phi(X^{t,x}_0)+ \int_0^s
f\left(r,X^{t,x}_r,\widehat{Y}^{t,x}_r \varepsilon_r(T_r),
\widehat{Z}^{t,x}_r\varepsilon_r(T_r)\right)\varepsilon_r^{-1}
(T_r)\d r-\int_0^s\widehat{Z}^{t,x}_r \downarrow \d W_r.$$ Let $0\le
t\le t'\le T,\, x,x'\in \R^d.$ Then by Proposition \ref{prop_meas},
we have, $$
\widehat{u}(t,x)-\widehat{u}(t',x')=E\left[\widehat{Y}^{t,x}_t-\widehat{Y}^{t',x'}_{t'}|
{\cal F}^B_T\right]$$ and, thus,
$$
\begin{aligned}
&\left|\widehat{u}(t,x)-\widehat{u}(t',x')\right|^p=
\left|E\left[\widehat{Y}^{t,x}_t-\widehat{Y}^{t',x'}_{t'}| {\cal
F}^B_T\right]\right|^p  \\
&\le
CE\left[\left|\widehat{Y}^{t,x}_t-\widehat{Y}^{t',x'}_{t}\right|^p
|{\cal F}^B_T\right]+
C\left|E\left[\widehat{Y}^{t',x'}_t-\widehat{Y}^{t',x'}_{t'}|{\cal
F}^B_T\right]\right|^p,
\end{aligned}
$$
where, $P$-a.s., by Lemma \ref{lemma_est},
$$
E\left[\left|\widehat{Y}^{t,x}_t-\widehat{Y}^{t',x'}_{t}\right|^p
|{\cal F}^B_T\right]\le C\exp\{pI_T^*\}
\left((1+|x|^p+|x'|^p)|t-t'|^{p/2}+|x-x'|^p\right),
$$
and $$
 \begin{aligned}
&\left|E\left[\widehat{Y}^{t',x'}_t-\widehat{Y}^{t',x'}_{t'}|{\cal
F}^B_T\right]\right|^p \\
\le&
\left|E\left[\int_t^{t'}\left|f\left(s,X^{t',x'}_s,\widehat{Y}^{t',x'}_s\varepsilon_s(T_s),
\widehat{Z}^{t',x'}_s\varepsilon_s(T_s)\right)\varepsilon_s^{-1}(T_s)\right|
\d s
|{\cal F}^B_T\right]\right|^p \\
\le&
\left|E\left[\int_t^{t'}C\left(\varepsilon_s^{-1}(T_s)f(s,0,0,0)
\left(1+\left|X_s^{t^\pr,x^\pr}\right|\right)+
\left|\widehat{Y}^{t',x'}_s\right|+
\left|\widehat{Z}^{t',x'}_s\right|\right)\d s|{\cal F}^B_T\right]\right|^p\\
\le&
C(1+|x'|^p)|t-t'|^p\exp\{pI_T^*\}+C|t-t'|^{p/2}\left(E\left[\int_t^{t'}\left(
\left|\widehat{Y}^{t',x'}_s\right|^2+\left|\widehat{Z}^{t',x'}_s\right|^2\right)\d
s|{\cal
F}^B_T\right]\right)^{p/2}\\
 \le&
C(1+|x'|^p)\exp\{pI_T^*\}|t-t'|^{p/2}.
\end{aligned}
$$
Consequently, for all $(t,x),(t',x')\in[0,T]\times \R^d,\, p\ge 1,$
$P$-a.s.,
$$|\widehat{u}(t,x)-\widehat{u}(t',x')|\le C\exp\{I_T^*\}
\left((1+|x|+|x'|)|t-t'|^{1/2}+|x-x'|\right).
$$
Hence,
$$
E\left[|\widehat{u}(t,x)-\widehat{u}(t',x')|^p\right]\le C_p
\left((1+|x|+|x'|)^p|t-t'|^{p/2}+|x-x'|^p\right),$$ and
Kolmogorov's continuity criterion gives the existence of a
continuous version of $\widehat{u}$.\hfill$\cajita$

\bigskip
Henceforth we denote by $\mathcal{L}$ the second-order differential
operator:
$$
\mathcal{L}: =\frac 12 {\rm{tr}}(\sigma\sigma^*(x)D^2_{xx})
+b(x)\nabla_x,
$$
and we consider  the following stochastic partial differential
equations:
\begin{equation}\label{spde_w}
\left\{ \begin{array}{ll} \d \hu(t,x)=\left[{\cal L}\hu(t,x)+
f\left(t,x,\hu(t,x)\ep_t(T_t),\nabla_x
\hu(t,x)\sigma(x)\ep_t(T_t)\right) \ep_t^{-1}(T_t)\right]\d t,& t\in[0,T];\\
\hu(0,x)=\Phi(x).&
\end{array}\right.
\end{equation}
and
\begin{equation}\label{spde_B}
\left\{ \begin{array}{ll} \d u(t,x) = \left[ \mathcal{L}u(t,x) +
f\left(t,x,u(t,x),\nabla_x u(t,x)\sigma(x)\right)\right]\d t
+\gamma_t  u(t,x)\d B_t,& t\in[0,T];\\
u(0,x)=\Phi(x).&
\end{array}\right.
\end{equation}
Our objective is to characterize  $\hu(t,x)=\hy^{t,x}_t$ and
$u(t,x)=\hu(A_t,t,x)\ep_t=Y_t^{t,x}$ as the  viscosity solutions of
the above stochastic partial differential equations (\ref{spde_w})
and (\ref{spde_B}), respectively.
\begin{nota}
In fact, equation $(\ref{spde_w})$ is  a partial differential
equation with random coefficients which can be solved pathwisely,
and the equation $(\ref{spde_B})$ is  a stochastic partial
differential equation driven by the fractional Brownian motion $B$.
\end{nota}

First we give the definition of a pathwise viscosity solution of
SPDE (\ref{spde_w}).
\begin{definicion} A real valued continuous
random field $\hu: \Omega^{\pr}\times[0,T]\times\R^d\mapsto\R$ is
called a pathwise viscosity solution  of equation $(\ref{spde_w})$
if there exists a subset $\overline{\Omega}^\pr$ of $\Omega^\pr$
with $P^\pr(\overline{\Omega}^\pr)=1$, such that for all
$\omega^\pr\in \overline{\Omega}^\pr$, $\hu(\omega^\pr,\cdot,\cdot)$
is a viscosity solution for the PDE $(\ref{spde_w})$ at
$\omega^\pr$.
\end{definicion}

For the definition of the viscosity solution, which is a well-known
concept by now, we refer to the \textit{User's Guide} by Crandall et
al. \cite{CIL}.

\bigskip
For the proof that $\hu(t,x)$ is the pathwise viscosity solution of
equation $(\ref{spde_w})$, we need the following two auxiliary
results.
\begin{lema}\label{comparison}
$($Comparison result.$)$ Let $(${\bf{H3}}$)$ hold. Let
$\left(\hy^1(\omega^\pr,\cdot),\hz^1(\omega^\pr,\cdot)\right)$ and
$\left(\hy^2(\omega^\pr,\cdot),\hz^2(\omega^\pr,\cdot)\right)$ be
the solutions of BSDE $(\ref{bsde})$ with coefficients $(\xi_1,
f_1)$ and $(\xi_2,f_2)$, respectively. Then, if $\xi_1\leq \xi_2$
and $f_1\leq f_2,$ it holds that $\hy^1_t(\omega^\pr,\cdot)\leq
\hy^2_t(\omega^\pr,\cdot)$, $t\in[0,T]$, $P^W$-a.s., for all
$\omega^\pr\in\widetilde{\Omega}^\pr$.
\end{lema}
\noindent\textit{Proof:} For the proof the reader is referred to the
comparison for BSDEs by Peng \cite{P1}, or also to Buckdahn and Ma
\cite{bu}. \hfill$\cajita$

\begin{lema}$($A priori estimate.$)$ Let
$\omega^\pr\in\widetilde{\Omega}^\pr$ and let
$\left(\hy^1,\hz^1\right)$ and $\left(\hy^2,\hz^2\right)$ be the
solutions of BSDE $(\ref{bsde})$  with coefficients $(\xi_1, f_1)$
and $(\xi_2,f_2)$, respectively, and put $\delta
\hy(\omega^\pr,\cdot)=\hy^1(\omega^\pr,\cdot)-\hy^2(\omega^\pr,\cdot)$,
$\delta \xi=\xi^1-\xi^2$ and  $\delta f_s(\omega^\pr,\cdot) =
(f_1-f_2) \left( s, \hy^2_s(\omega^\pr,\cdot) \ep_s(T_s,\omega^\pr),
\hz^2_s(\omega^\pr,\cdot)\ep_s(T_s,\omega^\pr) \right) \ep_s^{-1}
(T_s,\omega^\pr).$ Moreover, let $\widetilde {C}$ be the Lipschitz
constant of $f_1$. Then there exist $\beta \geq
\widetilde{C}(2+\nu^2)+\mu^2, \nu>0,\mu>0,$ such that $P^W$-a.s.,
for $0\leq s \leq T$, and for all
$\omega^\pr\in\widetilde{\Omega}^\pr$,
\begin{equation}\label{est_apriori}
E^W\left[\exp\{\beta(T-s)\}\left|\delta\hy_s(\omega^\pr,\cdot)
\right|^2\right]\leq E^W\left[\exp\{\beta
T\}\left|\delta\xi\right|^2+\frac{1}{\mu^2}\int^s_0
\exp\{\beta(T-r)\}
 |\delta f_r(\omega^\pr,\cdot)|^2\d r\right],
\end{equation}
\begin{equation}\label{est_deltay}
E^W\left[\sup_{0\leq r\leq
T}\left|\delta\hy_r(\omega^\pr,\cdot)\right|^2
+\int^T_0\left|\delta\hz_r(\omega^\pr,\cdot)\right|^2\d r\right]\leq
C E^W\left[|\delta\xi|^2+\int^T_0 |\delta f_r(\omega^\pr,\cdot)|^2\d
r\right].
\end{equation}
\end{lema}
\noindent\textit{Proof:} Let $\omega^\pr\in\widetilde{\Omega}^\pr$.
By applying It\^o's formula to $\left|\delta
\hy_s(\omega^\pr,\cdot)\right|^2 \exp\{\beta(T-s)\}$, we obtain
\[
\begin{aligned}
&\left|\delta\hy_s(\omega^\pr,\cdot)\right|^2\exp\{\beta(T-s)\}+\int^s_0\exp\{\beta(T-r)\}
\left(\beta\left|\delta\hy_r(\omega^\pr,\cdot)\right|^2+
\left|\delta\hz_r(\omega^\pr,\cdot)\right|^2\right)\d
r\\
&\quad \quad+2\int^s_0\exp\{\beta(T-r)\}
\delta\hy_r(\omega^\pr,\cdot)\delta\hz_r(\omega^\pr,\cdot)\dw \d W_r\\
=&|\delta\xi|^2\exp\{\beta T\} + 2 \int^s_0 \delta
\hy_r(\omega^\pr,\cdot)\exp\{\beta(T-r)\}
\Big[f_1\left(r,\hy^1_r(\omega^\pr,\cdot)\ep_r(T_r,\omega^\pr),
\hz^1_r(\omega^\pr,\cdot)\ep_r(T_r,\omega^\pr)\right) \\
&\quad \quad-
f_2\left(r,\hy^2_r(\omega^\pr,\cdot)\ep_r(T_r,\omega^\pr),
\hz^2_r(\omega^\pr,\cdot)\ep_r(T_r,\omega^\pr)\right) \Big]
\ep_r^{-1}(T_r,\omega^\pr) \d r,
\end{aligned}
\]
and by taking the expectation with respect to $P^W$ on both sides we
get, for $\nu>0, \mu>0$, $P$-a.s.,
\[
\begin{aligned}
&E^W\left[\exp\{\beta(T-s)\}\left|\delta\hy_s(\omega^\pr,\cdot)
\right|^2+\int^s_0\exp\{\beta(T-r)\}\left( \beta\left|
\delta\hy_r(\omega^\pr,\cdot)\right|^2+
\left|\delta\hz_r(\omega^\pr,\cdot)\right|^2\right)\d
r\right]\\
\leq&E^W\left[|\delta\xi|^2\exp\{\beta T\}\right]\\
&+E^W\left[ \int^s_0 \exp\{\beta(T-r)\} \left[ \big(
(2+\nu^2)\widetilde{C}+\mu^2\big)\left|\delta\hy_r(\omega^\pr,\cdot)\right|^2+
\frac{\widetilde{C} \left|\delta \hz_r(\omega^\pr,\cdot)\right|^2}
{\nu^2} +\frac{|\delta f_r(\omega^\pr,\cdot)|^2}{\mu^2}\right]\d r
\right].
\end{aligned}\]
Finally, by choosing $\beta\geq \widetilde{C}(2+\nu^2)+\mu^2$, with
$\nu^2>\widetilde{C}$, we obtain that $P$-a.s.,
$$
E^W\left[\exp\{\beta(T-s)\}\left|\delta\hy_s(\omega^\pr,\cdot)\right|^2\right]\leq
E^W\left[\exp\{\beta T\} \left|\delta\xi\right|^2+ \int^s_0
\exp\{\beta(T-r)\}\frac{1}{\mu^2}|\delta f_r(\omega^\pr,\cdot)|^2\d
r\right],
$$
which is exactly (\ref{est_apriori}).

Estimate (\ref{est_deltay}) can be proven by the arguments developed
for  (\ref{est:Yt-Yt'}).\hfill$\cajita$

\bigskip
Let us now turn to the solutions of our SPDEs. The next theorem is
one of the main results of this section. For this let
${\widehat{\Omega}^\pr}=\left\{\omega^{\pr}\in\widetilde{\Omega}^\pr|
\hu(\omega^\pr,\cdot,\cdot)\ \textrm{is continuous}\right\}$ and
notice that in light of Lemma \ref{continuity},
$P^B({\widehat{\Omega}}^\pr)=1$.
\begin{teorema}\label{thm-visc-w}
The random field  $\hu$ defined by
$\hu(\omega^{\pr},t,x)=\hy^{t,x}_t(\omega^\pr)$ for all
$\omega^\pr\in\widehat{\Omega}^\pr$ is a pathwise viscosity solution
of equation $(\ref{spde_w})$, where $\hy^{t,x}(\omega^\pr,\cdot)$ is
the solution of equation $(\ref{x-bsde})$. Furthermore, this
solution $\hu(\omega^{\pr},t,x)$ is unique in the class of
continuous stochastic fields $\tilde{u}: \Omega^\pr\times
[0,T]\times \R^d\mapsto \R$ such that, for some random variable
$\eta\in L^0(\mathcal{F}^B_T)$,
$$
|\tilde{u}(\omega^\pr,t,x)|\leq \eta(\omega^\pr)(1+|x|),  (t,x)\in
[0,T]\times\R^d, P^B(\d \omega^\pr)\textrm{-a.s.}
$$
\end{teorema}
\begin{nota}\label{remark_unique}
The uniqueness of the solution is to be understood as a $P$-almost
sure one: let $\hu_i (i=1,2)$ be such that $\hu_i(\omega^\pr)$ is a
viscosity solution of PDE $(\ref{spde_w})$ at $\omega^\pr$, for all
$\omega^\pr\in\widehat{\Omega}_i^\pr$. Then, by the uniqueness
result of viscosity solution of deterministic PDEs $($see: Pardoux
\cite{Pa}, Theorem 6.14$)$ we know that
$\hu_1(\omega^\pr,\cdot)=\hu_2(\omega^\pr,\cdot)$, for all
$\omega^\pr\in\widehat{\Omega}_1^\pr\cap\widehat{\Omega}_2^\pr$. In
particular, $\hu_1=\hu_2,$ $P$-a.s.
\end{nota}

\noindent\textit{Proof of Theorem \ref{thm-visc-w}:} We adapt the
method used in the paper of El Karoui et al.\cite{EPQ} to show that
$\hu$ is a pathwise viscosity subsolution of equation
(\ref{spde_w}). Recall that the set
$\widetilde{\Omega}^\pr:=\left\{\omega^\pr\in\Omega^\pr |
I_T^\ast(\omega^\pr)<+\infty\right\}$ satisfies
$P^B(\widetilde{\Omega}^\pr)=1$. We work here on the set
${\widehat{\Omega}^\pr}:=\big\{\omega^{\pr}\in\widetilde{\Omega}^\pr|
\hu(\omega^\pr,\cdot,\cdot)\ \textrm{is continuous}\big\}$ which
satisfies $P^B({\widehat{\Omega}}^\pr)=1$ in light of Lemma
\ref{continuity}.

Now, according to the definition of the viscosity solution, for an
arbitrarily chosen $\omega^\pr\in {\widehat{\Omega}^\pr}$, we fix
arbitrarily a point $(t,x)\in [0,T]\times \R^d$ and a test function
$\varphi\in C^{\infty}_b$ such that
$\varphi(t,x)=\hu\left(\omega^\pr,t,x\right)$ and $\varphi\geq
\hu\left(\omega^\pr,\cdot,\cdot\right)$.

For $t\in[0,T]$ and $h\geq0$, we have, thanks to the Proposition
\ref{prop_meas} and equation (\ref{x-bsde}),
$$
\hu\left(\omega^\pr,t,x\right)=\hu\left(\omega^\pr,t-h,X_{t-h}^{t,x}\right)+
\int^t _{t-h}
F_r\left(\omega^\pr,X_r^{t,x},\hy^{t,x}_r(\omega^\pr,\cdot), \hz^{t,
x}_r(\omega^\pr,\cdot)\right) \d r
-\int^{t}_{t-h}\hz^{t,x}_r(\omega^\pr,\cdot)\dw \d W_r.
$$
We emphasize that for fixed $\omega^{\pr}$, this BSDE can be viewed
as a classical BSDE with respect to $W$ and we recall that
$\hy^{t,x}_{t-h}\left(\omega^\pr,\cdot\right)
=\hu\left(\omega^\pr,t-h,X_{t-h}^{t,x}\right)$. Now for the fixed
$\omega^\pr\in {\widehat{\Omega}^\pr}$, it holds that
$$
\varphi(t,x)\leq\varphi\left(t-h,X_{t-h}^{t ,x}\right)+
\int^t_{t-h}F_r\left(\omega^\pr,X_r^{t,x},\hy^{t,x}_r(\omega^\pr,\cdot),
\hz^{t,x}_r(\omega^\pr,\cdot)\right)\d r
-\int^{t}_{t-h}\hz^{t,x}_r(\omega^\pr,\cdot)\dw \d W_r.
$$
Let $\left(\overline{Y}^{t,x,h}(\omega^\pr,\cdot),\overline{Z}^{t,
x,h}(\omega^\pr,\cdot)\right)\in
L^2_{\mathbb{G}}\left(t-h,t;\R\times\R^d\right)$ be the solution of
the following  equation evaluated at $\omega^{\pr}$: for
$s\in[t-h,t]$,
\begin{equation}\label{eq:overlineY}
\begin{aligned}
\overline{Y}^{t,x,h}_{s}(\omega^\pr,\cdot)=&\varphi\left(t-h,X_{t-h}^{t,
x}\right)+
\int^s_{t-h}F_r\left(\omega^\pr,X_r^{t,x},\overline{Y}^{t,x,h}_r(\omega^\pr,\cdot),
\overline{Z}^{t,x,h}_r(\omega^\pr,\cdot)\right)\d r\\
&-\int^{s}_{t-h}\overline{Z}^{t,x,h}_r(\omega^\pr,\cdot)\dw \d W_r.
\end{aligned}
\end{equation}
From Lemma \ref{comparison}, it follows that
$\overline{Y}^{t,x,h}_{t}(\omega^\pr,\cdot)\geq\varphi(t,x)
=\hu(\omega^\pr,t,x)$. Now we put
\begin{equation}\label{eq:wtYwtZ}
\wt{Y}_{s}^{t,x}(\omega^\pr,\cdot)=\overline{Y}_{s}^{t,x,h}(\omega^\pr,\cdot)-
\varphi\left(s,X_s^{t,x}\right)-\int^s_{t-h}G(\omega^\pr,r,x)\d r \
\end{equation}
and
\begin{equation*}
\wt{Z}_{s}^{t,x}(\omega^\pr,\cdot)=\overline{Z}_s^{t,x,h}(\omega^\pr,\cdot)-
(\nabla_x\varphi \sigma)(s,X_s^{t,x}),
\end{equation*}
where
$G(\omega^\pr,s,x)=\partial_s\varphi(s,x)-\mathcal{L}\varphi(s,x)-
F_s\big(\omega^\pr,x,\varphi(s,x),\nabla_x
\varphi(s,x)\sigma(x)\big).$ From the equations
(\ref{eq:overlineY}), (\ref{eq:wtYwtZ}) and It\^o's formula we have
\begin{equation*}
\begin{aligned}
\wt{Y}^{t,x}_{s}(&\omega^\pr,\cdot)= \int^s_{t-h}\bigg[
F_r\Big(\omega^\pr,X_r^{t,x},\varphi(r,X_r^{t,x}) +\wt{Y}^{t,
x}_r(\omega^\pr,\cdot)+\int^r_{t-h}G(\omega^\pr,s,x)\d s,
\\&(\nabla_x\varphi \sigma) (r,X_r^{t,x})+
\wt{Z}^{t,x}_r(\omega^\pr,\cdot)\Big)
-(\partial_r\varphi-\mathcal{L}\varphi)
(r,X_r^{t,x})+G(\omega^\pr,r,x)\bigg]\d r
-\int^{s}_{t-h}\wt{Z}^{t,x}_r(\omega^\pr,\cdot)\dw \d W_r.
\end{aligned}
\end{equation*}
Putting
$$\begin{aligned}
\delta(\omega^\pr,r,h)=&F_r\left(\omega^\pr,X_r^{t,x},\varphi(r,X_r^{t,x})
+\int^r_{t-h}G(\omega^\pr,s,x)\d s, (\nabla_x\varphi \sigma)
(r,X_r^{t,x})\right) \\
&-(\partial_r\varphi-\mathcal{L}\varphi)
(r,X_r^{t,x})+G(\omega^\pr,r,x),
\end{aligned}
 $$ we have
$|\delta(\omega^\pr,r,h)|\leq
\kappa(\omega^\pr)\left|X_r^{t,x}-x\right|, r \in [0,t],$ for some
$\mathcal{F}^B_T$-measurable and $P^B$-integrable $\kappa:
\Omega^\pr\mapsto\R^+$. From the a priori estimate
(\ref{est_apriori}), it follows that
\begin{equation}\label{est_delta}
\begin{aligned}
E^W\left[\sup_{t-h\leq s\leq t}\left|\wt{Y}_s^{t,
x}(\omega^\pr,\cdot)\right|^2+\int^t_{t-h}\left|\wt{Z}_r^{t,x}(\omega^\pr,\cdot)\right|^2\d
r\right]\leq C
E^W\left[\int^t_{t-h}\left|\delta(\omega^\pr,r,h)\right|^2\d
r\right]=h\rho(\omega^\pr,h).
\end{aligned}
\end{equation}
where $\rho(\omega^\pr,h)$  tends to 0 as $h\to 0$. Consequently, it
yields
\begin{equation}\label{est_visc_YZ}
E^W\left[\int^{t}_{t-h}\left(\left|\wt{Y}_r^{t,x}(\omega^\pr,\cdot)\right|+
\left|\wt{Z}_r^{t,x}(\omega^\pr,\cdot)\right|\right)\d
r\right]=h\sqrt{\rho}(\omega^\pr,h).
\end{equation}
Furthermore, we have
$\wt{Y}_{t}^{t,x}(\omega^\pr)=E^W\left[\wt{Y}_{t}^{t,x}(\omega^\pr,\cdot)\right]
=E^W\left[\int^t_{t-h}\delta^{\pr}(\omega^\pr,r,h)\d r\right]$,
 where
$$
\begin{aligned}
\delta^{\pr}(\omega^\pr,&r,h)=-(\partial_r\varphi-\mathcal{L}\varphi)
(r,X_r^{t,x})+G(\omega^\pr,r,x)\\&+F_r\left(\omega^\pr,X_r^{t,x},
\varphi(r,X_r^{t,x}) +\wt{Y}^{t, x}_r(\omega^\pr,\cdot)
+\int^r_{t-h}G(\omega^\pr,s,x)\d s, (\nabla_x\varphi
\sigma)(r,X_r^{t,x})+\wt{Z}^{t,x}_r(\omega^\pr,\cdot)\right).
\end{aligned}
$$
From the fact that $f$ is Lipschitz and the estimates
(\ref{est_delta}) and (\ref{est_visc_YZ}), we get
$$
\left|\wt{Y}_{t}^{t,x}(\omega^\pr)\right|\leq E^W\left[ \int^t_{t-h}
\left( |\delta(\omega^\pr,r,h)|
+\left|\wt{Y}_r^{t,x}(\omega^\pr,\cdot)\right|+
\left|\wt{Z}_r^{t,x}(\omega^\pr,\cdot)\right| \right)\d r
\right]=h\rho(\omega^\pr,h).
$$
Thus, from (\ref{eq:wtYwtZ}) (for $s=t$) and since
$\overline{Y}^{t,x}_t(\omega^\pr,\cdot)\geq \varphi(t,x)$, we obtain
$\int^{t}_{t-h}G(\omega^\pr,r,x)\d r\geq -h\rho(\omega^\pr,h)$.
Consequently, $\frac{1}{h}\int^{t}_{t-h}G(\omega^\pr,r,x)\d r\geq
-\rho(\omega^\pr,h)$. By letting $h$ tend to $0$, we finally get for
$\omega^\pr\in \widehat{\Omega}^\pr$,
$$
G(\omega^\pr,t,x)=\partial_t\varphi(t,x)- \mathcal{L} \varphi (t,x)
- F_{t}\big(\omega^\pr,x, \varphi(t,x),\nabla_x
\varphi(t,x)\sigma(x)\big)\geq0.
$$
Hence $\hu(\omega^\pr,t,x)$ is a  pathwise viscosity  subsolution of
$(\ref{spde_w})$. The proof of $\hu$ being a pathwise viscosity
supersolution is similar.

The proof of uniqueness becomes clear from Remark
\ref{remark_unique}.\hfill$\cajita$

\bigskip
In analogy to the relation between the solutions $(Y,Z)$ and
$(\hy,\hz)$ of the associated BDSDE and the BSDE, respectively, we
shall expect that $u(t,x)=Y^{t,x}_t=\hy^{t,x}_t(A_t)\ep_t,
(t,x)\in[0,T]\times \R^d$, is a solution of SPDE (\ref{spde_B}).
This claim is confirmed by
\begin{proposicion}\label{thm-spde-B}
Suppose that $u, \hu$ are $C^{0,2}$-stochastic fields over
$\Omega^\pr\times[0,T]\times\R^d$ such that there exist $\delta>0$
and a constant $C_{\delta,x}>0$ $($only depending on $\delta$ and
$x$$)$ with:
\begin{equation}\label{condition}
E\left[|w(t,x)|^{2+\delta}+\int^t_0\left(|\nabla_x
w(s,x)|^{2+\delta} +|D^2_{xx}w(s,x)|^{2+\delta}\right)\d
s\right]\leq C_{\delta,x}, t\in[0,T],\ \textrm{for} \ w=u,\hu.
\end{equation}

 Then $\hu(t,x)$ is a  classical pathwise
solution of equation $({\ref{spde_w}})$ if and only if $u(t,x)$ is a
classical solution of SPDE $(\ref{spde_B})$.
\end{proposicion}
\noindent\textit{Proof:}  We restrict ourselves to show $u(t,x)$
solves equation $(\ref{spde_B})$ whenever $\hu$ solves
(\ref{spde_w}). For this, we proceed in analogy to the proof of
Theorem \ref{thm-BDSDE}. Let $F$ be an arbitrary but fixed element
of $\mathcal{S}_{\mathcal{K}}^{\pr}$. By using the Girsanov
transformation we have
\begin{equation*}
\begin{aligned}
&E\left[u(t,x)F-u(0,x)F\right]= E\left[ \hu(A_t,t,x) \ep_t F -
\hu(0,x) F\right] =E\left[F(T_t)\hu(t,x)-F\hu(0,x)\right]\\
=&E\left[F(T_t)\hu(0,x)-F\hu(0,x)\right]\\
&+E\left[F(T_t)\int^t_0\left[{\cal L}\hu(s,x)+ f \big( s, x,
\hu(s,x) \ep_s(T_s), \nabla_x\hu(s,x)\sigma(x)\ep_s(T_s)\big)
\ep_s^{-1} (T_s) \right]\d s\right].
\end{aligned}
\end{equation*}
As in the proof of Theorem \ref{thm-BDSDE} we use the fact that
$\frac{\d}{\d t} F(T_t)=\gamma_t (\mathcal{K}^*\mathcal{K}D^B
F)(t,T_t)$ to deduce the following
\begin{equation*}
\begin{aligned}
&E\left[u(t,x)F-u(0,x)F\right]\\=&E\left[\hu(0,x)
\int^t_0\gamma_s(\mathcal{K}^{\ast}\mathcal{K}D^BF(T_s))(s)\d
s\right]+E\left[\int^t_0\int^t_s\gamma_r(\mathcal{K}^{\ast}
\mathcal{K}D^BF(T_r))(r)\d r{\cal L}\hu(s,x)\d s\right]\\
&+E\left[\int^t_0F(T_s)\left[{\cal L}\hu(s,x)+
f\big(s,x,\hu(s,x)\ep_s(T_s),\nabla_x \hu(s,x) \sigma(x)
\ep_s(T_s) \big)\ep_s^{-1}(T_s)\right]\d s\right]\\
&+ E\left[\int^t_0\int^t_s\gamma_r(\mathcal{K}^{\ast}
\mathcal{K}D^BF(T_r))(r)\d rf\big(s,x,\hu(s,x)\ep_s(T_s),\nabla_x
\hu(s,x)\sigma(x)\ep_s(T_s)\big) \ep_s^{-1}(T_s)\d s\right].
\end{aligned}
\end{equation*}
Thanks to the assumption that $\hu(t,x)$ is a pathwise classical
solution of equation $({\ref{spde_w}})$, we obtain
\begin{equation*}
\begin{aligned}
&E\left[u(t,x)F-u(0,x)F\right]\\=&E\left[\int^t_0F(T_s)\left[{\cal
L}\hu(s,x)+ f\big(s,x,\hu(s,x)\ep_s(T_s),\nabla_x \hu(s,x) \sigma(x)
\ep_s(T_s)\big)
\ep_s^{-1}(T_s)\right]\d s\right]\\
&+E\left[\int^t_0\gamma_s(\mathcal{K}^{\ast}\mathcal{K}D^BF(T_s))(s)
\hu(s,x)\d s\right]\\
=&E\left[\int^t_0F\left[{\cal L}{u}(s,x)+ f \big( s, x, {u}(s,x),
\nabla_x {u}(s,x)\sigma(x)\big) \right]\d s\right]+ E \left[
\int^t_0 \gamma_s(\mathcal{K}^{\ast}\mathcal{K}D^BF)(s) {u}(s,x)\d
s\right].
\end{aligned}
\end{equation*}
Consequently,
\begin{equation*}\
\begin{aligned}
&E\left[\int^t_0\gamma_s(\mathcal{K}^{\ast}\mathcal{K}D^BF)(s)
{u}(s,x)\d s\right]
\\=&E\left[F\Big(u(t,x)-u(0,x)-\int^t_0\left[{\cal L}{u}(s,x)+ f
\big( s, x, {u}(s,x),\nabla_x {u}(s,x)\sigma(x)\big)\right]\d s
\Big) \right].
\end{aligned}
\end{equation*}
From the integrability condition (\ref{condition}) we know that
$$
u(t,x)-u(0,x)-\int^t_0\left[{\cal L}{u}(s,x)+ f \big( s, x,
{u}(s,x), \nabla_x {u}(s,x)\sigma(x)\big)\right]\d s\in L^2 (\Omega,
\mathcal{F},P).
$$
Moreover, $\gamma1_{[0,t]}u\in L^2([0,T]\times\Omega)$. Indeed,
$$E\left[\int^T_0|\gamma_r1_{[0,t]}(r)u(r,x)|^2\d r\right]=
\int^T_0|\gamma_r1_{[0,t]}(r)|^2E[|u(r,x)|^2]\d r\leq
C_{\delta,x}\int^T_0|\gamma_r1_{[0,t]}(r)|^2\d r<\infty.$$ By
Definition ${\ref{def-div-ext}}$ we get
\begin{equation*}
\begin{aligned}
&E\left[F\int^t_0\gamma_s{u}(s,x)\d B_s\right]
=E\left[F\Big(u(t,x)-u(0,x)-\int^t_0\left[{\cal L}{u}(s,x)+ f \big(
s, x,{u}(s,x),\nabla_x {u}(s,x)\sigma(x)\big)\right]\d s \Big)
\right].
\end{aligned}
\end{equation*}
It then follows from the arbitrariness of $F\in
\mathcal{S}^{\pr}_{\mathcal{K}}$ that
\begin{equation*}
u(t,x)=\Phi(x)+\int^t_0\left[{\cal L}{u}(s,x)+
f\big(s,x,{u}(s,x),\nabla_x {u}(s,x)\sigma(x)\big)\right] \d
s+\int^t_0\gamma_s {u}(s,x)\d B_s.
\end{equation*}
The proof is complete now.\hfill$\cajita$
\begin{Nota}
The regularity  of $\hu$ in the above proposition is difficult to
get under not too restrictive assumptions $($like coefficients of
class $C^3_{l,b}$, linearity of $f$ in $z$$)$.
\end{Nota}
\begin{nota}
Notice that generally speaking, a continuous random field after
Girsanov  transformation $A_t$ is not necessarily continuous in $t$
any more.
We give a simple counterexample:\\
 Let $0<s<T$ be fixed and
$$
\hu(t,x)= \left\{
\begin{array}{ll}
(t-s)B_t,& 0\leq t \leq s;\\
(t-s)\emph{sgn}(B_s),& s< t\leq T.
\end{array}\right.
$$
It is obvious that $\hu(t,x)$ is $\mathcal{F}^B_t$-measurable and
continuous in $t$. But after Girsanov transformation $A_t$, it
becomes
$$
u(t,x)=\hu(A_t,t,x)\ep_t=\left\{
\begin{array}{ll}
(t-s)\left(B_t-\int^t_0(\mathcal{K}\gamma\textbf{1}_{[0,t]})(r)\d
r\right)\ep_t,& 0\leq t \leq s;\\
(t-s)\emph{sgn}\left(B_s - \int^s_0 (\mathcal{K} \gamma \textbf{1}
_{[0,t]}) (r)\d r \right)\ep_t,& s< t\leq T,
\end{array}\right.
$$
which is not continuous in $t$ on
$$\left\{\omega^\pr:\inf_{t\in[s,T]}\left(B_s-\int^s_0\left(\K \gamma
1_{[0,t]}\right)(r)\d
r\right)<0<\sup_{t\in[s,T]}\left(B_s-\int^s_0\left(\K \gamma
1_{[0,t]}\right)(r)\d r\right)\right\}.$$
\end{nota}
However, as we state below, the random field $u$  has a continuous
version in our case. To this end we need the following technical
result:
\begin{lema}\label{lemma_conti}
Let $\gamma$ be such that $({\bf H1})$ holds. Then there exist
positive constants $C$ and $q$ such that for all
$r,v,v^\pr\in[0,T]$, $v\le v^\pr$,
\[
\left|\int^T_0({\cal K}\gamma 1_{[v^\pr,v]})(s)({\cal K}\gamma
1_{[0,r]})(s)\d s\right|\leq C |v-v^\pr|^q.
\]
\end{lema}
\noindent\textit{Proof:} We have
\[
\begin{aligned}
&\left|\int^T_0({\cal K}\gamma 1_{[v^\pr,v]})(s)({\cal K}\gamma
1_{[0,r]})(s)\d s\right| \leq \left(\int^T_0(\K \gamma
1_{[v^\pr,v]})^2(s)\d s\right)^{1/2}\left(\int^T_0(\K \gamma
1_{[0,r]})^2(s)\d s\right)^{1/2}\\ \leq &C\left(\int^T_0(\K \gamma
1_{[v^\pr,v]})^2(s)\d s\right)^{1/2},
\end{aligned}
\]
where the last inequality follows from \cite{Ls} (Lemma 2.3). Also,
by the proof of Lemma 2.3 in \cite{Ls} and using the notation
$\alpha=1/2-H$, we have
\begin{equation}\label{eq:I1+I2}
\begin{aligned}
&(\K \gamma 1_{[v^\pr,v]})(s) \\=&1_{[v^\pr,v]}(s)
\left(\phi_{\gamma}(s)+\frac{\alpha
s^{\alpha}}{\Gamma(1-\alpha)}\int^T_v\frac{r^{-\alpha}\gamma_r}
{(r-s)^{1+\alpha}}\d r\right)-1_{[0,v^\pr]}(s)\frac{\alpha
s^{\alpha}}{\Gamma(1-\alpha)}\int^v_{v^\pr}\frac{r^{-\alpha}\gamma_r}
{(r-s)^{1+\alpha}}\d r\\
=&I_1(s)+I_2(s).
\end{aligned}
\end{equation}
Now, applying \cite{Ls} (Lemma 2.3) again, we obtain
\begin{equation}
\begin{aligned}
&\int^T_0I_1(s)^2\d s =\int^T_01_{[v^\pr,v]}(s)
\left(1_{[0,v]}(s)\left[\phi_{\gamma}(s)+\frac{\alpha
s^{\alpha}}{\Gamma(1-\alpha)}\int^T_v\frac{r^{-\alpha}\gamma_r}
{(r-s)^{1+\alpha}}\d r\right]\right)^2\d s \\
\leq &
|v-v^\pr|^{(p^\pr-2)/p^\pr}\|\phi_{\gamma1_{[0,v]}}\|^2_{L^{p^\pr}}
\leq C|v-v^\pr|^{(p^\pr-2)/p^\pr}.
\end{aligned}
\end{equation}

On the other hand, for $m=1+\eta$ and $q_m=m/\eta,$ with $\eta$
small enough, we can write
\[
\begin{aligned}
&\int^v_{v^\pr}\frac{r^{-\alpha}\gamma_r} {(r-s)^{1+\alpha}}\d r\\
\leq&|v-v^\pr|^{{1}/{q_m}}\left[\int^v_{v^{\pr}}
\frac{r^{-m\alpha}|\gamma_r|^m}{(r-s)^{m(1+\alpha)}}\d
r\right]^{1/m}\\
=
&\frac{1}{\Gamma(\alpha)}|v-v^\pr|^{{1}/{q_m}}\left[\int^v_{v^{\pr}}
\frac{\left(\int^T_r\phi_{\gamma}(\theta)\theta^{-\alpha}
(\theta-r)^{\alpha-1}\d \theta\right)^m}{(r-s)^{m(1+\alpha)}}\d
r\right]^{1/m}\\
\leq &C|v-v^\pr|^{{1}/{q_m}}\left[\int^T_{v^{\pr}}
|\phi_{\gamma}(\theta)|^m\theta^{-\eta}
\int^{\theta}_{v^\pr}(r-s)^{-m(1+\alpha)}r^{-m\alpha+\eta}
(\theta-r)^{m(\alpha-1)} \d r
\d \theta\right]^{1/m}\\
\leq
&C|v-v^\pr|^{{1}/{q_m}}(v^\pr-s)^{-{2\eta}/{m}}\left[\int^T_{v^{\pr}}
|\phi_{\gamma}(\theta)|^m\theta^{-\eta}
\int^{\theta}_{v^\pr}(r-s)^{-m\alpha+\eta-1}r^{-m\alpha+\eta}
(\theta-r)^{m\alpha-\eta-1} \d r \d \theta\right]^{1/m}.
\end{aligned}
\]
Hence, \cite{Ls} (Lemma 2.2) gives
\[\begin{aligned}
&\int^v_{v^\pr}\frac{r^{-\alpha}\gamma_r} {(r-s)^{1+\alpha}}\d r\\
\leq&
C|v-v^\pr|^{{1}/{q_m}}(v^\pr-s)^{-{2\eta}/{m}}\left[\int^T_{v^{\pr}}
|\phi_{\gamma}(\theta)|^m\theta^{-\eta}
{v^\pr}^{-m\alpha+\eta}(v^\pr-s)^{-m\alpha+\eta}(\theta-s)^{-1}
(\theta-v^\pr)^{m\alpha-\eta} \d \theta\right]^{{1}/{m}}\\
\leq& C|v-v^\pr|^{{1}/{q_m}}
(v^\pr-s)^{-\alpha-{\eta}/{m}}\left[\int^T_{v^{\pr}}
|\phi_{\gamma}(\theta)|^m\theta^{-\eta}
{v^\pr}^{-m\alpha+\eta}(\theta-s)^{-1+m\alpha-\eta}\d
\theta\right]^{1/m},
\end{aligned}\]
which, together with (\ref{eq:I1+I2}), implies, for $p<1/\alpha,$ $
p^\pr=\frac{pm}{2(1-(m\alpha-\eta)p)}$ and
${1}/{p^\pr}+{1}/{q^\pr}=1$,
\[\begin{aligned}
&\int^T_0I_2(s)^2\d s\\
\leq &C |v-v^\pr|^{\frac{2}{q_m}}
\left(\int^{v^\pr}_0(v^\pr-s)^{2q^\pr(-\alpha-{\eta}/{m})}\d s
\right)^{\frac{1}{q^\pr}}\left(\int_0^{v^\pr}\left[\int^T_{v^{\pr}}
|\phi_{\gamma}(\theta)|^m\theta^{-\eta}
(\theta-s)^{-1+m\alpha-\eta}\d \theta\right]^{\frac{2p^\pr}{m}}\d
s\right)^{\frac{1}{p^\pr}}\\
\leq& C|v-v^\pr|^{{2}/{q_m}}.
\end{aligned}
\]
Thus, we get the wished result.\hfill$\cajita$

\begin{lema}\label{lemma_conti_u}
The random field $u(t,x) :=\hu(A_t,t,x)\ep_t, (t,x) \in [0, T]
\times \R^d$ has a continuous version.
\end{lema}
\noindent\textit{Proof:} In the following, for simplicity of
notations, we put
$$\Theta_r^{t,x,v,v^\pr}=\left(r,X_r^{t,x},\hy_r^{t,x}(A_v) \ep_r(T_rA_{v^\pr}),
\hz_r^{t,x}(A_v) \ep_r(T_rA_{v^\pr})\right).$$ For $0\leq v^\pr\leq
v\leq T$, we notice that $\left(\hy^{t,x}(A_v)\right)$ is the
solution of the BSDE
\begin{equation}\label{eq:Y_Av}
\begin{aligned}
\hy_s^{t,x}(A_v)= \Phi\left(X_0^{t,x}\right) +\int^s_0
f\left(\Theta^{t,x,v,v}_r\right)\ep_r^{-1}(T_rA_v)\d r -\int_0^s
\hz_r^{t,x}(A_v)\dw \d W_r,
\end{aligned}
\end{equation}
while $(\hy^{t,x}(A_{v^\pr}))$ is the solution of the BSDE
\[
\begin{aligned}
\hy_s^{t,x}(A_{v^\pr})=\Phi\left(X_0^{t,x}\right) +\int^s_0
f\left(\Theta^{t,x,v^{\pr},v^{\pr}}_r\right)\ep_r^{-1}(T_rA_{v^\pr})\d
r -\int_0^s \hz_r^{t,x}(A_{v^\pr})\dw \d W_r.
\end{aligned}
\]
We set $J_r^v=\exp\left\{\int^T_0({\cal K}\gamma 1_{[0,v]})(s)({\cal
K}\gamma 1_{[0,r]})(s)\d s\right\}$, and we observe that $
\ep^{-1}_r (T_rA_v)= \ep_r^{-1}(T_r)J^v_r$. Moreover, equation
(\ref{eq:Y_Av}) can be written  as follows:
\begin{equation*}
\begin{aligned}
\hy_s^{t,x}(A_v)= &\Phi\left(X_0^{t,x}\right) +\int^s_0
f\left(r,X_r^{t,x},\left(\hy_r^{t,x}(A_v),
\hz^{t,x}_r(A_v)\right)\ep_r(T_r)(J_r^v)^{-1}\right)\ep_r^{-1}(T_r)J_r^v\d
r \\
&-\int_0^s \hz_r^{t,x}(A_v)\dw \d W_r,\quad s\in[0,t],
\end{aligned}
\end{equation*}
and a comparison with (\ref{x-bsde}) suggests the similarity of
arguments which can be applied. So, by a standard BSDE estimate,
see, for instance, the proof of Lemma \ref{lemma_est}, we have
\begin{equation}\label{est_deltay_Au}
\begin{aligned}
&E\left[\sup_{0\leq s \leq T}
\left|\hy^{t,x}_s(A_v)-\hy^{t,x}_s(A_{v^\pr})\right|^p\right] \\
\leq &CE\left[
\int^{T}_0\left|f(\Theta^{t,x,v,v}_r)\ep^{-1}_r(T_rA_v)
-f(\Theta^{t,x,v,v^{\pr}}_r)\ep^{-1}_r(T_rA_{v^\pr})\right|^p\d r\right]\\
\leq&CE\bigg[
\int^T_0\bigg[\left|f(\Theta^{t,x,v,v}_r)\left(\ep^{-1}_r(T_rA_v)-
\ep^{-1}_r(T_rA_{v^{\pr}})\right)\right|^p
+\left|\left(f(\Theta^{t,x,v,v}_r)
-f(\Theta^{t,x,v,v^{\pr}}_r)\right)\ep^{-1}_r(T_rA_{v^\pr})\right|^p\bigg]\d
r\bigg].
\end{aligned}
\end{equation}
Recalling that $ \ep^{-1}_r (T_rA_v)= \ep_r^{-1}(T_r)J^v_r$ and
applying Lemma \ref{lemma_conti}, we get
\begin{equation}\label{est_delta_ep}
\begin{aligned}
&\left|\ep_r^{-1}(T_rA_v)-\ep^{-1}_r(T_rA_{v^\pr})\right|^p=
\ep_r^{-p}(T_r)\left|J_r^{v}- J_r^{v^\pr}\right|^p\\
\leq& C\ep_r^{-p}(T_r)\left|\int^T_0({\cal K}\gamma
1_{[0,v]})(s)({\cal K}\gamma 1_{[0,r]})(s)\d s-\int^T_0({\cal
K}\gamma 1_{[0,v^\pr]})(s)({\cal K}\gamma 1_{[0,r]})(s)\d
s\right|^p \\
=&C\ep_r^{-p}(T_r)\left|\int^T_0({\cal K}\gamma
1_{[v^\pr,v]})(s)({\cal K}\gamma 1_{[0,r]})(s)\d s\right|^p\leq C
\ep_r^{-p}(T_r) \left|v-v^\pr\right|^{pq}.
\end{aligned}
\end{equation}
On the other hand,
\begin{equation}\label{est_deltaf}
\begin{aligned}
&\Big|f(\Theta^{t,x,v,v}_r)
-f(\Theta^{t,x,v,v^{\pr}}_r)\Big|^p\ep^{-p}_r(T_rA_{v^\pr})\\
\leq &C\left(\left|\hy^{t,x}_r(A_v)\right|^p+\left|\hz^{t,x}_r(A_v
)\right|^p\right)
\left|\ep_r(T_rA_v)-\ep_r(T_rA_{v^{\pr}})\right|^p\ep_r^{-p}
(T_rA_{v^{\pr}})\\
=&C\left(\left|\hy^{t,x}_r(A_v)\right|^p+\left|\hz^{t,x}_r(A_v)
\right|^p\right) \left|\exp\left\{-\int^T_0({\cal K}\gamma
1_{[v^\pr,v]})(s)({\cal K}\gamma 1_{[0,r]})(s)\d s\right\}-1\right|^p\\
\leq& C\left(\left|\hy^{t,x}_r(A_v)\right|^p+\left|\hz^{t,x}_r(A_v)
\right|^p\right) \left|v-v^\pr\right|^{pq}.
\end{aligned}
\end{equation}
Plugging estimates (\ref{est_delta_ep}) and (\ref{est_deltaf}) into
the equation (\ref{est_deltay_Au}), Lemma \ref{lemma_est} yields
that
\[
\begin{aligned}
&E\left[\sup_{0\leq s \leq
T}\left|\hy^{t,x}_s(A_v)-\hy^{t,x}_s(A_{v^\pr})\right|^p\right] \\
\leq&CE\bigg[
\int^T_0\left(\left|f(\Theta^{t,x,v,v}_r)\left(\ep^{-1}_r(T_r)
\right)\right|^p+\left|\hy^{t,x}_r(A_v)\right|^p+\left|\hz^{t,x}_r(A_v)
\right|^p\right)\d r\bigg] \left|v-v^\pr\right|^{pq}\\
\leq&C(1+|x|^p)\left|v-v^\pr\right|^{pq}.
\end{aligned}
\]
For the latter inequality we have used the following estimate of
$\hz$, which proof will be postponed until the end of the current
proof.
\begin{lema}\label{lemma_z} There exists a constant $C_p$ such that
$E\left[\sup_{0\leq s \leq T} \left|\hz^{t,x}_s\right|^p\right]\leq
C_p(1+|x|^p).$
\end{lema}

We continue our proof of Lemma \ref{lemma_conti_u}:  According to
the proof of Lemma \ref{continuity} we have
\[
\begin{aligned}
&E\left[\sup_{0\leq s\leq T}
\left|\hy^{t^\pr,x^\pr}_s(A_v)-\hy^{t,x}_s(A_v)\right|^p \right]
=E\left[E\left[\sup_{s\in[0,T]}\left|\hy^{t^\pr,x^\pr}_s-\hy^{t,x}_s
\right|^p\big|
\mathcal{F}_T^B\right]\circ A_v\right]\\
\leq&
\left(E\left[\sup_{s\in[0,T]}\left|\hy^{t^\pr,x^\pr}_s-\hy^{t,x}_s\right|^{2p}
\right]\right)^{1/2}
\left(E\left[\ep_v^{-2}(T_v)\right]\right)^{1/2} \leq
C\left((1+|x|^p+|x^\pr|^p)|t-t^\pr|^{p/2}+|x-x^\pr|^p\right).
\end{aligned}
\]
Hence, by combining the above estimates we obtain
\[
\begin{aligned}
&E\left[\sup_{0\leq s\leq T}\left|\hy^{t^\pr,x^\pr}_s(A_{v^\pr}) -
\hy^{t,x}_s(A_v)\right|^p\right] \\
\leq& C\left(E\left[\sup_{0\leq s\leq
T}\left|\hy^{t^\pr,x^\pr}_s(A_{v^\pr}) -
\hy^{t^\pr,x^\pr}_s(A_{v})\right|^p\right]+ E\left[\sup_{0\leq s\leq
T}\left|\hy^{t^\pr,x^\pr}_s(A_v) -
\hy^{t,x}_s(A_v)\right|^p\right]\right)\\
\leq
&C\left((1+|x|^p+|x^\pr|^p)(|t-t^\pr|^{p/2}+|v-v^\pr|^{pq})+|x-x^\pr|^p\right).
\end{aligned}
\]
Consequently, from the Kolmogorov continuity criterion we know the
process $\{\hy^{t,x}_s(A_v);s,t,v\in[0,T],x\in\R^d\}$ has an a.s.
continuous version. From Lemma \ref{lemma_conti} we have that
$\ep_t$ is continuous in $t$. It then follows that
$u(t,x)=\hy^{t,x}_t(A_t)\ep_t$ has a version which is jointly
continuous in $t$ and $x$.\hfill$\cajita$

\bigskip

\noindent\textit{Proof of Lemma \ref{lemma_z}:} For simplicity we
suppose all the functions $\Phi, f, \sigma, b$ are smooth. For the
proof of the general case, the Lipschitz functions have to be
approximated by smooth functions with the same Lipschitz constants.
We define $(\underline{Y}^{t,x}, \underline{Z}^{t,x})$ to be the
solution of the equation
\begin{equation}\label{eq:L^2 dy}
\underline{Y}^{t,x}_s=\Phi^{\pr}(X_0^{t,x})\nabla_x X_0^{t,x}+
\int^s_0 \left[f^{\pr}_x(\Xi_r^{\pr})\ep_r^{-1}(T_r) \nabla_x
X_r^{t,x}+f^{\pr}_y(\Xi_r^{\pr})\underline{Y}^{t,x}_r +f_z^{\pr}(
\Xi_r^{\pr})\underline{Z}_r^{t,x}\right]\d r
-\int^s_0\underline{Z}_r^{t,x}\dw \d W_r,
\end{equation}
where $\Xi_r^{\pr}=\left(r,X_r^{t,x},\hy_r^{t,x}\ep_r(T_r),
\hz_r^{t,x}\ep_r(T_r)\right)$ and
\begin{equation}\label{eq:nabla-x}
\nabla_x X_r^{t,x}=I - \int^t_r \nabla_x
X_s^{t,x}b^\pr\left(X_s^{t,x}\right)\d s -\int^t_r \nabla_x
X_s^{t,x}\sigma^\pr\left(X_s^{t,x}\right)\dw \d W_s, \ r\in[0,t].
\end{equation}
With the arguments used in the proof of Lemma 4.2, we get
\begin{equation}\label{est_Y_}
E\left[\sup_{0\leq s\leq T}\left| \underline{Y}^{t,x}_s\right|^p +
\left( \int^T_0|\underline{Z}^{t,x}_s|^2\d s\right)^{p/2}\right]\leq
C_p(1+|x|^p).
\end{equation}
By arguments which by now are standard, it can be seen that the
processes $X^{t,x}, \hy^{t,x}$ and $\hz^{t,x}$ are Malliavin
differentiable with respect to $W$, and thus, for $s\leq u\leq t\leq
T$,
\begin{equation}\label{eq_D^W Y}
\begin{aligned}
D^W_u\hy^{t,x}_s=&\Phi^{\pr}(X_0^{t,x})D_u^W X_0^{t,x}
-\int^s_0D_u^W\hz_r^{t,x}\dw \d W_r\\&+\int^s_0
\left[f_x^{\pr}(\Xi_r^{\pr})\ep_r^{-1}(T_r) D_u^W
X_r^{t,x}+f^{\pr}_y(\Xi_r^{\pr})D_u^W\hy^{t,x}_r +f_z^{\pr}(
\Xi_r^{\pr})D^W_u\hz_r^{t,x}\right]\d r.
\end{aligned}
\end{equation}
On the other hand, from
$$
\hy^{t,x}_s=\hy^{t,x}_{\theta}+\int^s_{\theta}f\left(r,X_r^{t,x},
\hy_r^{t,x}\ep_r(T_r), \hz_r^{t,x}\ep_r(T_r)\right)\ep_r^{-1}(T_r)\d
r-\int^s_{\theta}\hz_{r}^{t,x}\dw\d W_r, \theta< s\leq u\leq T,
$$
we get
\begin{equation}\label{eq:Z}
\begin{aligned}
D^W_u\hy^{t,x}_s=&-\hz^{t,x}_u +\int^u_sD_u^W\hz_r^{t,x}\dw \d
W_r\\&+\int^u_s \left[f_x^{\pr}(\Xi_r^{\pr})\ep_r^{-1}(T_r) D_u^W
X_r^{t,x}+f^{\pr}_y(\Xi_r^{\pr})D_u^W\hy^{t,x}_r +f_z^{\pr}(
\Xi_r^{\pr})D^W_u\hz_r^{t,x}\right]\d r,
\end{aligned}
\end{equation}
From the above three equations (\ref{eq:L^2 dy}), (\ref{eq_D^W Y})
and (\ref{eq:Z}) and the relation
\begin{equation}\label{eq:DX}
D^W_r(X_s^{t,x})= \nabla_x X_s^{t,x}(\nabla_x
X_r^{t,x})^{-1}\sigma(X_r^{t,x})
\end{equation}
(see: Pardoux and Peng \cite{PP2}), we obtain that, for $s\leq u\leq
t\leq T$,
$$
D^W_u \hy_s^{t,x}=\underline{Y}_s^{t,x}(\nabla_x
X_u^{t,x})^{-1}\sigma(X_u^{t,x}),
$$
$$
\hz_s^{t,x}=-\underline{Y}_s^{t,x}(\nabla_x
X_s^{t,x})^{-1}\sigma(X_s^{t,x}).
$$
Finally, from standard estimates for (\ref{eq:nabla-x}) and from the
estimate (\ref{est_Y_}), we get
$$\begin{aligned}
E\left[\sup_{0\leq s \leq T} \left|\hz^{t,x}_s\right|^p\right]\le&
E\left[\sup_{0\le s\le T} \left|\underline{Y}^{t,x}_s\right|^p
\sup_{0\le s\le T}\left|\nabla_x X_s^{t,x}\right|^{-p} \sup _{0\le
s\le
T}\left|\sigma(X_s^{t,x})\right|^p\right]\\
\le&E\left[\sup_{0\le s\le T}
\left|\underline{Y}^{t,x}_s\right|^{pq_1}\right]^{{1}/{q_1}}
E\left[\sup_{0\le s\le T}\left|\nabla_x
X_s^{t,x}\right|^{-pq_2}\right]^{{1}/{q_2}}
E\left[\sup _{0\le s\le T}\left|\sigma(X_s^{t,x})\right|^{pq_3}\right]^{{1}/{q_3}}\\
\leq& C_p(1+|x|^p),
\end{aligned}
$$
where ${1}/{q_1}+{1}/{q_2}+{1}/{q_3}=1,$ with $q_1, q_2, q_3>1$. The
proof is complete.\hfill$\cajita$

\bigskip
The above proposition motivates the following definition:
\begin{definicion}
A  continuous random field $u: [0,T] \times \R^d \times \Omega^{\pr}
\mapsto\R$ is  a $($stochastic$)$ viscosity solution of equation
$(\ref{spde_B})$ if and only if $\hu(t,x)=u(T_t,t,x)\ep^{-1}_t(T_t),
(t,x)\in[0,T]\times\R^d$ is a pathwise viscosity  solution of
equation $(\ref{spde_w})$.
\end{definicion}
As a consequence of our preceding discussion, we can formulate the
following statement:
\begin{teorema}
The continuous stochastic field
$u(t,x):=\hu(A_t,t,x)\ep_t=\hy^{t,x}_t(A_t)\ep_t=Y^{t,x}_t$ is a
stochastic viscosity solution of the semilinear SPDE
$(\ref{spde_B})$. This solution is unique inside the class $C_p^B$
of continuous stochastic field $\tilde{u}: \Omega^\pr\times
[0,T]\times \R^d\mapsto \R$ such that,
$$|\tilde{u}(t,x)|\leq
C\exp\{I_T^\ast\}(1+|x|), (t,x)\in[0,T]\times \R^d, P-a.s.,
$$
for some constant $C$ only depending on $\tilde{u}$.
\end{teorema}
\begin{nota}
$1)$ It can be easily checked that $\tilde{u}\in C_p^B$ if and only
if $(\tilde{u}(A_t,t,x)\ep_t)\in C_p^B$ if and only if $
(\tilde{u}(T_t,t,x)\ep^{-1}_t(T_t))\in C_p^B
$.\\
$2)$ As a consequence of the preceding theorem we have that
$u(t,x)=Y^{t,x}_t$ is the unique $($in $C_p^B$$)$ stochastic
viscosity solution of SPDE $(\ref{spde_B})$. This extends the
Feynman-Kac formula to SPDEs driven by a fractional Brownian motion.
\end{nota}

We conclude the main theorems of Section \ref{sec:4} with the
following relation diagram, which shows the mutual relationship
between fractional backward SDEs and SPDEs:
$$
\begin{array}{|ccc|}
\hline \hu(t,x)\ is\ the\ viscosity\ solution\ of\
(\ref{spde_w})&\stackrel{GT}{\longleftrightarrow}&u(t,x)\ is\ the\
viscosity \
solution\ of\ (\ref{spde_B})\\
\updownarrow& &\updownarrow\\
(\hy,\hz)\ is\ the\ solution\ of\ BSDE\
(\ref{x-bsde})&\stackrel{GT}{\longleftrightarrow}&(Y,Z)\ is\  the\
solution\ of\ BDSDE\ (\ref{x-bdsde})\\
\hline
\end{array}
$$
where $'GT'$ stands for $'$Girsanov transformation$'$.

\bigskip
Finally, in order to illustrate how our method works, we give the
example of a linear fractional backward doubly stochastic
differential equation.
\begin{ejemplo}
We let $d=1$ and $f(s,x,y,z)=f^1_s x+f^2_s y+f^3_sz$, where the
coefficients $f^1_s, f^2_s$ and $f^3_s$ are bounded and
deterministic functions. The associated fractional backward doubly
stochastic differential equation is linear and writes:
\begin{equation}\label{eq:example}
Y_s^{t,x}=\Phi(X_0^{t,x})+\int^s_0 \left(f^1_r X_r^{t,x}+f^2_r
Y^{t,x}_r+f^3_r Z^{t,x}_r\right)\d r -\int^s_0 Z^{t,x}_r\dw \d W_r
+\int^s_0 \gamma_r Y_r^{t,x}\d B_r.
\end{equation}
After Girsanov transformation, it becomes
\[
\hy_s^{t,x}=\Phi(X_0^{t,x})+\int^s_0 \left(f^1_r
X_r^{t,x}\ep_r^{-1}(T_r)+f^2_r \hy^{t,x}_r+f^3_r
\hz^{t,x}_r\right)\d r -\int^s_0 \hz^{t,x}_r\dw \d W_r.
\]
and has the following solution:
\[
\hy^{t,x}_s=E_Q\left[\int^s_0\left(f_r^1X_r^{t,x}\ep_r^{-1}(T_r)
\exp\left\{\int^s_rf^2_u\d u\right\}\right)\d r+
\exp\left\{\int^s_0f_r^2\d
r\right\}\Phi(X_0^{t,x})\Big|\mathcal{F}_{s,t}^{W}
\vee\mathcal{F}_t^B\right],
\]
where $E_Q$ is the expectation with respect to
$Q=\exp\left\{\int^t_0 f_r^3\d W_r-\frac12\int^t_0(f_r^3)^2\d
r\right\}P.$ According to Theorem {\ref{thm-BDSDE}}, the solution of
$(\ref{eq:example})$ is then
\[
Y^{t,x}_s=E_Q\left[\ep_s\int^s_0\left(f_r^1X_r^{t,x}\ep_r^{-1}(T_rA_s)
\exp\left\{\int^s_rf^2_u\d u\right\}\right)\d r+ \ep_s \exp \left\{
\int^s_0 f_r^2\d r\right\} \Phi (X_0^{t,x}) \Big| \mathcal{F}_{s,t}
^W \vee\mathcal{F}_t^B\right].
\]
\end{ejemplo}
\section*{Acknowledgements} Part of this work was done while
the second author was visiting  the Institut Mittag-Leffler. He
thanks it for its hospitality and economical support. The authors
also thank Professor Rainer Buckdahn for his valuable advice and
discussions.

\end{document}